\documentclass[journal,onecolumn,12pt]{IEEEtran}
\usepackage[numbers,sort,compress]{natbib}
\usepackage{pslatex}
\usepackage{amsfonts,color,morefloats}
\usepackage{amssymb,amsmath,latexsym,amsthm}
\usepackage{hyperref}
\usepackage{color}
\usepackage{makecell}
\usepackage{bm}

\usepackage[utf8]{inputenc}
\usepackage[T1]{fontenc}  

\DeclareUnicodeCharacter{202F}{\,} 

\hypersetup{hidelinks}
\allowdisplaybreaks

\newcommand{\gf}{{\mathbb{F}}}
\newcommand{\F}{{\mathbb{F}}}

\newcommand{\bH}{{\mathrm{H}}}

\newtheorem{theorem}{Theorem}
\newtheorem{lemma}{Lemma}
\newtheorem{proposition}{Proposition}
\newtheorem{corollary}{Corollary}

\newtheorem{remark}{Remark}

\begin{document}
\title{On covering radius of generalized Zetterberg codes}
\author{Maosheng Xiong\IEEEauthorrefmark{1}, Haode Yan\IEEEauthorrefmark{2} \\
\IEEEauthorblockA{\IEEEauthorrefmark{1}Department of Mathematics, The Hong Kong University of Science and Technology, Hong Kong, P. R. China. \\
\IEEEauthorrefmark{2}School of Science, Harbin Institute of Technology, Shenzhen, 518055, China. \\ 
E-mail: \href{mailto: mamsxiong@ust.hk}{mamsxiong}@ust.hk, \href{mailto: hdyan@swjtu.edu.cn}{hdyan}@swjtu.edu.cn}\\
	}
\maketitle
\begin{abstract} 
We employ analytic number-theoretic techniques---specifically, character sums and Weil‑type estimates---to study the covering radius of the generalized Zetterberg codes $C_s(q_0)$ of length $q_0^{s}+1$ and dimension $\,q_0^{s}+1-2s\,$ over the finite field $\mathbb{F}_{q_0}$. Although the even- and odd-field cases require distinct technical treatment, the proofs follow a unified analytic framework that is substantially simpler and more transparent than previous approaches. We prove that $\rho(C_s(q_0)) \le 3$ for all $q_0$ and $s$, and determine its exact value for a wide range of $s$ for each $q_0$. When $q_0$ is even, our results fill the gap left by recent studies focused solely on odd $q_0$; for odd $q_0$, the range of $s$ for which $\rho(C_s(q_0))$ is exactly determined is considerably broader than previously known. Combined with the corresponding minimum distance results, we obtain infinitely many quasi‑perfect and maximal codes within this family.
\end{abstract}
	
	

\thispagestyle{empty}   

\section{Introduction} \label{sec1}

\subsection{The covering radius}

Let $\mathbb{F}_{q_0}$ denote the finite field of $q_0$ elements, where $q_0$ is a prime power. Let $n$ be a positive integer. For two vectors $\bm{x}=(x_1,\cdots,x_n), \bm{y}=(y_1,\cdots,y_n)\in \mathbb{F}_{q_0}^n$, their \emph{Hamming distance} is 
\[
d_\bH(\bm{x}, \bm{y}) = \#\{\, i : x_i \neq y_i \,\}. 
\]
Here $\#I$ denotes the cardinality of a finite set $I$. Any subset 
$C \subseteq \mathbb{F}_{q_0}^n$ is called a \emph{code} of length $n$ over $\mathbb{F}_{q_0}$. The \emph{covering radius} of $C$ is defined as 
\[
\rho(C) = \max_{\bm{x} \in \mathbb{F}_{q_0}^n} d_\bH(\bm{x}, C), \quad d_\bH(\bm{x}, C) = \min_{\bm{c} \in C} d_\bH(\bm{x}, \bm{c}).\] 
Equivalently, $\rho(C)$ is the smallest integer $r$ such that  the Hamming spheres of radius $r$ centered at the codewords cover $\F_{q_0}^n$. This parameter provides a fundamental geometric measure of how densely the codewords of $C$ occupy the ambient space.

The covering radius is a fundamental concept in coding theory and discrete mathematics. 
It arises naturally in \textbf{decoding}, where~$\rho(C)$ represents the maximum number of errors a received word may contain before becoming ambiguous under nearest-neighbor decoding; 
in \textbf{data compression} and \textbf{quantization}, where it corresponds to the worst-case distortion under nearest-codeword approximation and relates to sphere-covering problems; 
and in \textbf{testing}, \textbf{fault detection}, and \textbf{memory coding}, where it measures the completeness of test sets and the robustness of storage schemes.
Comprehensive treatments of these topics can be found in the standard monographs and surveys~\cite{BLP98,CHLL97,CKMS85,CLLM97}. 
In \textbf{finite geometry}, the covering radius is closely related to the study of complete caps in projective spaces; 
see, for instance,~\cite{MR1714367,MR4568358,MR3809962,MR3156928} for representative research contributions and further developments.

Beyond these applications, the covering radius lies at the core of several notions of optimality in coding theory, giving rise to perfect, quasi‑perfect, and maximal codes---classes that represent extremal configurations of the Hamming space. The following discussion recalls these definitions and their relation to $\rho(C)$.

It is well known that
\[
\left\lfloor \frac{d(C)-1}{2} \right\rfloor \le \rho(C). 
\]
A code satisfying 
\[\rho(C)=\left\lfloor \frac{d(C)-1}{2} \right\rfloor\]
is called \emph{perfect}; such codes partition the Hamming space into disjoint Hamming spheres and simultaneously achieve optimal error-correction capability and maximum packing efficiency. When 
\begin{eqnarray} \label{x1:qp} \rho(C)=\left\lfloor \frac{d(C)-1}{2} \right\rfloor+1,\end{eqnarray}
the code is termed \emph{quasi-perfect}, extending the perfect packing by one additional layer. 

The parameters of all perfect codes over finite fields were completely determined in the 1970s~\cite{T73,T74}, and the classification of linear perfect codes is now fully known (see \cite{HP03}). Only a few exceptional families---notably the Hamming and Golay codes---attain this equality. By constrast, the classification of quasi-perfect codes remains substantially more complicated.  Consequently, recent research has focused on constructing explicit families of quasi-perfect codes with diverse parameters~(see, e.g., \cite{Danev,DDR11,D86,D85,1522651,GPZ60,H78,LH16}); see also~\cite{SHOS22,SHO23,SLHO25} for several of the latest developments.

A code $C$ is \emph{maximal} if it cannot be strictly contained in a larger code with the same minimum distance.  Equivalently \cite{CKMS85}, a maximal code satisfies 
\begin{eqnarray} \label{x1:mc}
\rho(C) \le d(C)-1,
\end{eqnarray}
meaning that every vector of the Hamming space lies within distance $d(C)-1$ from $C$. In this sense, the covering radius links two opposite extremal behaviors: perfect codes realize optimal packing, filling the space with disjoint spheres, while maximal codes represent the most dense configurations possible under the given distance constraint. Thus $\rho(C)$ quantifies how closely a code approaches these two geometric extremes within $\mathbb{F}_{q_0}^n$.

Despite extensive study, determining the covering radius of a general linear code remains a challenging problem---it is both \textbf{NP-hard} and  \textbf{co-NP-hard}~\cite{M84}. Consequently, most work has focused on deriving upper and lower bounds rather than exact values~\cite{AB02,H85,MC03,S86,MR3592808,MR3876436,MR3768641,MR2092631,MR1431812,MR994108}).  Exact covering radii are known only for a few highly structured families of codes; see, for example, \cite{Danev,DDR11,D86,D85,GPZ60,H78,H96,LH16}. For most recent progress, we refer the reader to  \cite{MR4604461,MR4709188,MR4395429,SLHO25}. The generalized Zetterberg codes introduced  below provide another notable instance in which precise analysis is attainable.

\subsection{Generalized Zetterberg Codes}

Among the many families of linear codes studied for their distance and covering properties, \emph{Zetterberg codes} occupy a distinctive place. Introduced independently by  Zetterberg~\cite{Z62} and
Meggitt~\cite{meggittErrorCorrecting1961}, these cyclic codes combine a high information rate with double-error-correcting capability and
admit efficient decoding algorithms~\cite{dodunekovAlgebraicDecoding1992}.
They are quasi-perfect~\cite{D85}, reversible~\cite[Ch.~7,
Sec.~6]{MS77}, and their $\mathbb{Z}_4$-lifts preserve these algebraic features~\cite{alahmadiLiftedZetterberg2016}. Owing to their balance of theoretical elegance and practical efficiency, Zetterberg codes remain an active research  topic~\cite{jingResultZetterberg2010,DHNX}.

Originally defined over $\mathbb{F}_2$, these codes extend naturally to any finite field~$\mathbb{F}_{q_0}$, giving rise to the \emph{generalized Zetterberg codes}~$C_s(q_0)$. 
They retain the same cyclic structure while revealing rich number-theoretical connections: for instance, the weight distributions of binary and ternary Zetterberg codes 
relate to Hecke operators acting on certain spaces of some modular forms
\cite{SCHOOF1991163,VANDERGEER1991256}, and the number of weight-five codewords connects to the Hasse zeta function of certain K3~surfaces~\cite{petersHasseZeta1992}. Such links place  $C_s(q_0)$ at the intersection of coding theory, number theory, and arithmetic geometry, and they also appear in the study of Niho-type cross-correlation \cite{xiaOpenProblem2016,xiaCorrelationDistribution2017,
	xiongCorrelationDistribution2024a,ShuxingXiongYan}.


\subsection{Summary of Main Results}

Recent research~\cite{SHO23,SLHO25} has provided a detailed account of the covering radius of generalized Zetterberg codes~$C_s(q_0)$ for odd~$q_0$. Although these works established impressive results, their proofs rely on lengthy and intricate elementary arguments, and the \emph{even-field case}~($q_0$ even) has remained largely open. 

In this paper we employ analytic number-theoretic techniques---specifically, character sums and Weil-type estimates---to study the covering radius of $C_s(q_0)$ for both even and odd $q_0$. This approach greatly advances the understanding of the even-field case and, for the odd-field case, yields new and stronger results than previously obtained. The main results are summarized below: 
\begin{enumerate}
\item[(1)] We prove that the covering radius of $C_s(q_0)$ satisfies $\rho(C_s(q_0)) \le 3$ for all $q_0$ and $s$;
\item[(2)] We establish necessary and sufficient conditions under which $\rho(C_s(q_0))=3$, and determine its exact value for a wide range of $s$ for each $q_0$; when $q_0$ is odd, this range extends substantially beyond that of \cite{SHO23,SLHO25}; 
	\item[(3)] By analyzing the parameters of $C_s(q_0)$, we identify infinitely many quasi-perfect and maximal codes arising from this construction. 
\end{enumerate}

Although the proofs for even‑ and odd-field cases require different arguments, the overall analytic framework remains uniform. Compared with \cite{SHO23,SLHO25}, our proofs are simpler and more transparent, leading to stronger results for odd~$q_0$ and substantial new progress on the even-field case, though a broad range of $s$ for which $\rho(C_s(q_0))$ is undetermined remains open for future investigation.

This paper is organized as follows. Section \ref{Notation} introduces the notation and basic conventions. In Section \ref{sec2} we define the generalized Zetterberg code $C_s(q_0)$ and the half-Zetterberg code $\widetilde{C_{s}(q_0)}$ for odd $q_0$, and collect several key lemmas used later in the analysis. Section \ref{sec3} investigates the minimum distance of $C_s(q_0)$ for even $q_0$ and that of $\widetilde{C_{s}(q_0)}$ for odd $q_0$. In Section \ref{sec4} we prove that $\rho(C_s(q_0)) \le 3$ for all $q_0$ and $s$, using character sums and Weil-type estimates; the arguments for even and odd $q_0$ require separate treatments. A detailed study of $\rho(C_s(q_0))$ for odd $q_0$ is given in Section \ref{sec5}, followed by the case of even $q_0$ in Section \ref{sec6}. Section \ref{sec7} combines our results on the minimum distance and covering radius to identify many families of quasi-perfect and maximal codes derived from $C_s(q_0)$. Finally, in Section \ref{sec8} we conclude the paper and propose an open problem.

\section{Notation} \label{Notation}

We use the following notation throughout the paper.
 
 \begin{itemize}
 	\item $q_0$ denotes a prime power, $s$ a positive integer, $q=q_0^s$; 
 	\item $\gf_{q_0},\gf_{q}$ and $\gf_{q^2}$ are the finite fields of order $q_0,q$ and $q^2$ respectively;
 	\item $\gf_q^*:=\gf_q \setminus \{0\}$;
 	\item If $q$ is odd, $\square_q \subseteq \gf_q^*$ denotes the set of squares;  
 	\item $H \subseteq \gf_{q^2}$ is the multiplicative subgroup of order $q+1$, defined by
 	\[H:=\left\{x \in \gf_{q^2}: x^{q+1}=1\right\};\]
 	\item For any $x \in \gf_{q^2}$, write $\overline{x}:=x^q$;
 	\item $C_s(q_0)$ denotes the $q_0$-ary generalized Zetterberg code defined in \eqref{x:zett}; 
 	\item For a code $C$, let $d(C)$ and $\rho(C)$ denote its minimum distance covering radius respectively; 
 	\item for a positive integer $n$, $[n]:=\left\{1,2,\cdots,n\right\}$;
 	\item $\mathbb{P}^1(\gf_q)=\gf_q \cup \{\infty\}$. 
 \end{itemize}

\section{Preliminaries}\label{sec2}

\subsection{Generalized Zetterberg and half-Zetterberg codes}

Let $\xi$ be a generator of $H$ so that 
\[H=\left\{1,\xi,\xi^2,\cdots, \xi^{q}\right\}.\]
The $q_0$-ary \emph{generalized Zetterberg code}, denoted by $C_s(q_0)$, is defined as 
\begin{eqnarray} \label{x:zett} C_s(q_0)=\left\{(c_0,c_1,\cdots,c_q) \in \gf_{q_0}^{q+1}: \sum_{i=0}^q c_i \xi^i=0\right\}.\end{eqnarray}
Thus $C_s(q_0)$ is a $q_0$-ary linear code of length $q+1$. 

Let $h_{q_0}(X) \in \gf_{q_0}[X]$ be the monic irreducible polynomial of $\xi$ over $\gf_{q_0}$. Since $\gf_{q^2}=\gf_q(\xi)$, we have
\[\left[\gf_q(\xi):\gf_{q_0}\right]=\left[\gf_q(\xi):\gf_{q}\right]\cdot \left[\gf_q:\gf_{q_0}\right]=2s,\]
so $\deg h_{q_0}(X)=2s$ and $h_{q_0}(X)|\left(x^{q+1}-1\right)$.

Each $\underline{c}=(c_0,\ldots,c_q)\in \gf_{q_0}^{q+1}$ corresponds to a polynomial $\underline{c}(X)=c_0+c_1X+\cdots+c_qX^q$. This gives an isomorphism from the vector space $\gf_{q_0}^{q+1}$ to the quotient ring $\gf_{q_0}[X]/\left(X^{q+1}-1\right)$. It is easy to see that  
\[\underline{c} \in C_s(q_0) \Longleftrightarrow \underline{c}(\xi)=0 \Longleftrightarrow h_{q_0}(X)|\underline{c}(X),\]
hence $C_s(q_0)$ can be identified with the cyclic group in $\gf_{q_0}[X]/\left(X^{q+1}-1\right)$ generated by $h_{q_0}(X)$. Since $\deg h_{q_0}(X)=2s$, we have 
\[\dim C_s(q_0)=q+1-2s,\]
so the code $C_s(q_0)$ has parameters $\left[q+1,q+1-2s\right]_{q_0}$. 

We will carry out a detailed study of the minimum distance of $C_s(q_0)$ for $q_0$ even in Section \ref{sec3}. When $q_0$ is odd, the situation is much easier: since $q+1$ is even, $\pm 1 \in H$, so the minimum distance of $C_s(q_0)$ is always 2. The minimum distance can be improved by dropping half the elements of $H$. To be more precise, let us define the $q_0$-ary \emph{half-Zetterberg code} $\widetilde{C_s(q_0)}$ (when $q_0$ is odd) as 
\begin{eqnarray} \label{x:halfZ} \widetilde{C_s(q_0)}=\left\{\left(c_0,c_1,\cdots,c_{{(q-1)}/{2}}\right) \in \gf_{q_0}^{{(q+1)}/{2}}: \sum_{i=0}^{{(q-1)}/{2}} c_i \xi^i=0\right\}.\end{eqnarray}
By the correspondence 
\[\underline{c}=\left(c_0,\ldots,c_{(q-1)/2}\right)\in \gf_{q_0}^{(q+1)/2} \mapsto \underline{c}(X)=c_0+c_1X+\cdots+c_{(q-1)/2}X^{(q_0-1)/2},\]
it is easy to see that 
\[\underline{c} \in \widetilde{C_s(q_0)} \Longleftrightarrow \underline{c}(\xi)=0 \Longleftrightarrow h_{q_0}(X)|\underline{c}(X),\]
hence $\widetilde{C_s(q_0)}$ can be identified with the constacyclic code in $\gf_{q_0}[X]/\left(X^{(q+1)/2}+1\right)$ generated by $h_{q_0}(X)$. Since $\deg h_{q_0}(X)=2s$, the code $\widetilde{C_s(q_0)}$ has parameters $\left[\frac{q+1}{2},\frac{q+1}{2}-2s\right]_{q_0}$ with minimum distance at least 3. 

We remark that in \cite{SHO23} the authors defined the so-called ``half and twisted Zetterberg codes'', but they are all monomially equivalent to $\widetilde{C_s(q_0)}$. Actually, for any $\underline{\epsilon}=\left(\epsilon_0,\cdots,\epsilon_{{(q-1)}/{2}}\right) \in \{\pm 1\}^{{(q+1)}/{2}}$, denote $\xi_i':=\epsilon_i \xi^i$ for $0 \le i \le {(q-1)}/{2}$ and define
\[\widetilde{C_s(q_0)}^{\underline{\epsilon}}:=\left\{\left(c_0,c_1,\cdots,c_{{(q-1)}/{2}}\right) \in \gf_{q_0}^{{(q+1)}/{2}}: \sum_{i=0}^{{(q-1)}/{2}} c_i \xi_i'=0\right\}.\]
It is easy to see that $\phi: \widetilde{C_s(q_0)} \to \widetilde{C_s(q_0)}^{\underline{\epsilon}}$ given by 
\[\left(c_0,c_1, \cdots,c_{{(q-1)}/{2}}\right)\mapsto \left(\epsilon_0c_0,\epsilon_1c_1, \cdots,\epsilon_{\frac{q-1}{2}}\,c_{\frac{q-1}{2}}\right)\]
is an isomorphism, hence the codes $\widetilde{C_s(q_0)}$ and $\widetilde{C_s(q_0)}^{\underline{\epsilon}}$ for any $\underline{\epsilon} \in \{\pm 1\}^{\frac{q+1}{2}}$ have the same parameters such as dimension, minimum distance and covering radius. It suffices to consider only the half-Zetterberg code $\widetilde{C_s(q_0)}$ given in \eqref{x:halfZ} when $q_0$ is odd. 

By the code parameters of $C_s(q_0)$ and $\widetilde{C_s(q_0)}$ and the sphere-packing bound (see \cite{HP03, MS77}) 
    \[A_{q_0}(n,d) \le \frac{q_0^n}{\sum_{i=0}^{\lfloor\frac{d-1}{2}\rfloor}\binom{n}{i}(q_0-1)^i},\]
where $A_{q_0}(n,d)$ denotes the largest cardinality of a $q_0$-ary code of length $n$ with minimum distance $d$,  we can easily obtain the following: 
\begin{eqnarray} \label{x:spbound} \left\{\begin{array}{lll}
d\left(C_s(q_0)\right) \le 4&:& q_0 \ge 4 \mbox{ even};\\
d\left(\widetilde{C_s(q_0)}\right) \le 4&:& q_0 \ge 5 \mbox{ odd}. 
\end{array} \right.\label{x:sphere} \end{eqnarray}
We will determine $d\left(C_s(q_0)\right)$ for $q_0$ even and $d\left(\widetilde{C_s(q_0)}\right)$ for $q_0$ odd in Section \ref{sec3}. 

As for the covering radius, it is worth mentioning that (see \cite{SHO23})
\[\rho\left(C_s(q_0)\right)=\rho\left(\widetilde{C_s(q_0)}\right), \quad q_0 \mbox{ is odd,} \]
and $\rho\left(C_s(q_0)\right)$ can be described in a very simple way (see \cite{CHLL97,H85} and \cite{SHO23}):
\begin{lemma} \label{x:crlemma} The covering radius of $C_s(q_0)$ is the least positive integer $\rho$ such that for any $\alpha \in \gf_{q^2}$, there exist $\left(c_1,\ldots,c_{\rho}\right) \in \gf_{q_0}^{\rho}$ and $\left(x_1,\ldots,x_{\rho}\right) \in H^{\rho}$ such that 
\[\alpha=c_1x_1+\cdots+c_{\rho} x_{\rho}. \]
\end{lemma}

\subsection{Character sums over finite fields}

Let $m,n$ be positive integers such that $m|n$. The trace and the norm functions from $\gf_{q_0^n}$ to $\gf_{q_0^m}$ are define as
\begin{align*}\mathrm{Tr}_{\gf_{q_0^n}/\gf_{q_0^m}}(x)&=x+x^{q_0^m}+\cdots+x^{q_0^{m\left(\frac{n}{m}-1\right)}},\\
\mathrm{N}_{\gf_{q_0^n}/\gf_{q_0^m}}(x)&=x \cdot x^{q_0^m} \cdots x^{q_0^{m\left(\frac{n}{m}-1\right)}}. \end{align*}

If $q$ is odd, let $\chi(\cdot)$ be the quadratic character of $\gf_q$, which is defined as
$$\chi(x)=
\left\{
\begin{array}{cl}
	1, & \mbox{if } x \in \square_q ,\\
	-1,& \mbox{if } x \in \gf_q^* \setminus \square_q, \\
	0, & \mbox{if } x=0.
\end{array}
\right.
$$
Let $\gf_q[X]$ be the polynomial ring over $\gf_q$. For $f\in\gf_q[X]$, we consider character sums of the form
\begin{eqnarray}\label{chasum} \sum_{x\in\gf_q}\chi(f(x)).\end{eqnarray}
The case that $\deg(f) =1$ is trivial. For $\deg(f) =2$, we have the following result:

\begin{lemma} \label{x:lem-qua} \cite[Theorem 5.48]{LN97} Let $f(X)=a_2X^2+a_1X+a_0 \in \gf_q[X]$ with $q$ odd and $a_2 \ne 0$. Put $d=a_1^2-4a_0a_2$ and let $\chi$ be the quadratic character of $\gf_q$. Then
\[\sum_{x \in \gf_q} \chi\left(f(x)\right)=\left\{
\begin{array}{cl}
    -\chi(a_2), & \mbox{ if } d \ne 0,\\
    (q-1) \chi(a_2), & \mbox{ if } d=0. 
\end{array}\right.\]
\end{lemma}

For a polynomial $f$ with degree $3$ or higher, computing the character sum (\ref{chasum}) is generally challenging. The subsequent lemma provides a standard Weil-type bound which is sufficient for our purpose. 
\begin{lemma}\cite[Corollary 2.3]{MR1401947}\label{weil}
	Let $f_1(T), \ldots,f_n(T)$ be $n$ monic pairwise prime polynomials in $\gf_q[T]$ whose largest squarefree divisors have degrees $d_1,\ldots,d_n$ respectively. Let $\chi_1,\ldots,\chi_n$ be multiplicative non-trivial characters of order $r_1,\ldots,r_n$ respectively on the finite field $\gf_q$. Assume that for some $1 \le i \le n$, the polynomial $f_i(T)$ is not of the form $g(T)^{r_i}$ in $\gf_q[T]$. Then, we have the estimate
	\begin{eqnarray} \label{x:wanlem} \left|\sum_{c\in\gf_q}\chi_1\left(f_1(c)\right)\cdots \chi_n\left(f_n(c)\right)\right|\leqslant\left(\sum_{i=1}^n d_i-1\right)\sqrt{q}.\end{eqnarray}
    Moreover, if $\chi_i^{d_i}=1$ for all $i$, then the right hand side of \eqref{x:wanlem} can be improved to 
    \[1+\left(\sum_{i=1}^n d_i-2\right)\sqrt{q}.\]
\end{lemma}

As a simple application of Lemma \ref{weil}, we obtain the following result.
\begin{lemma}\label{exist} Let $q_0\geq 5$ be an odd prime power. Then there exist $(c_1,c_2) \in \gf_{q_0}^{*2}$, such that
	\[\chi((c_1+c_2+1)(c_1+c_2-1)(c_1-c_2+1)(c_1-c_2-1))=-1,\]
	where $\chi(\cdot)$ denotes the quadratic multiplicative character of $\gf_{q_0}$.
\end{lemma}

\begin{proof} For an arbitrary $c_2\in\gf_{q_0}^*$, define the quartic polynomial  $f(X)=(X+c_2+1)(X+c_2-1)(X-c_2+1)(X-c_2-1) \in \gf_{q_0}[X]$. For $i=-1,0,1$, let
	\[n_i=\left|\left\{x\in\gf_{q_0}|\chi(f(x))=i\right\}\right|.\]
Clearly, $n_{-1}+n_0+n_{1}=q_0$, and since $f(X)$ has at most four roots in $\gf_{q_0}$, we have $0 \leq  n_0\leq 4$. Moreover,  because $c_2\neq 0$, $f$ is not a square of a polynomial over $\gf_{q_0}$. By Lemma \ref{weil}, we obtain
		$$\left|\sum_{x\in\gf_{q_0}}\chi\left(f(x)\right)\right|\leq 3\sqrt{q_0}.$$
Note that 
$$\sum_{x\in\gf_{q_0}}\chi\left(f(x)\right)=n_1-n_{-1},$$
this implies
\[\left|n_1-n_{-1}\right|\leq 3\sqrt{q_0}.\]
Combining this with $n_1+n_{-1}=q_0-n_0\geq q_0-4$, we deduce that 
\[n_{-1}\geq \frac{1}{2}\left(q_0-3\sqrt{q_0}-4\right).\]

Now for $q_0\geq 17$, a direct computation shows that $\frac{1}{2}\left(q_0-3\sqrt{q_0}-4\right)>0$, this with $\chi(f(0))\neq -1$ implies that there must exist some $c_1\in\gf_{q_0}^*$ such that $\chi(f(c_1))=-1$. For the remaining cases $q_0\in\{5,7,11,13\}$, the existence of such $c_1$ and $c_2$ in $\gf_{q_0}^*$ can be verified easily by using {\sc Magma}. This completes the proof of Lemma \ref{exist}. 
\end{proof}

An additive character of $\gf_{q}$ is a function $\psi: \gf_{q}\to\mathbb{C}\setminus \{0\}$ such that 
\[\psi(\beta_1+\beta_2)=\psi(\beta_1)\psi(\beta_2), \quad \forall \beta_1,\beta_2 \in \gf_q.\]
Suppose $\gf_q$ is of characteristic $p$. The canonical additive character is given by
\[\psi(\beta)=\zeta_p^{\mathrm{Tr}_{\gf_q/\gf_p}\left(\beta\right)}, \quad \forall \beta \in \gf_q, \]
where $\zeta_p:=\exp\left(2 \pi \sqrt{-1}/p\right)$ is the complex primitive $p$-th root of unity. When $p=2$, then $\zeta_2=-1$, this character takes values in $\{\pm 1\}$, a fact frequently used in coding theory. The following Weil-type bound was given in \cite{CP06}. 
\begin{lemma}\cite{CP06}\label{rational}
	 Consider the finite field $\gf_q$ with characteristic $p$. Let $f(X)\in \gf_q(X)$ be a rational function, 
	\[f(X)=p(X)+\frac{r(X)}{q(X)}, \quad \deg(r(X))<\deg(q(X)), \quad M:=\deg(p(X)),\]
	where $p(X),q(X),r(X)\in \gf_q[x]$ are polynomials, and $\gcd(r(X),q(X))=1$. Suppose further that $f(X)$ is non-constant, and that $M=0$ or $p\nmid M$. Write
	\[q(X)=\prod_{i=1}^Qq_i(X)^{m_i},\]
	where $q_i(X)$ are distinct irreducible polynomials over $\gf_q$, and $m_i\geq 1$ for each $i$. Define 
	\[L:=\sum_{i=1}^Q(m_i+1)\deg(q(X)).\]
	Then there exist complex numbers $\omega_j\in\mathbb{C}$, $1\leq j \leq M+L-1$, such that
	\[\sum_{\beta\in\gf_q\setminus S}\psi(f(\beta))=-\sum_{j=1}^{M+L-1}\omega_j,\]
	where $S$ is the set of poles of $f(X)$. Additionally, $\left|\omega_j\right|=\sqrt{q}$ for all $j$, except for a single value, $j'$, satisfying $|\omega_{j'}|=1$. Thus, 
	\[\left|\sum_{\beta\in\gf_q\setminus S}\psi(f(\beta))\right| \leq 1+\left(M+L-2\right)\sqrt{q}.\]
\end{lemma}

We will also use the following three criteria for quadratic polynomial over $\gf_q$. 
\begin{lemma} \label{x:chan-X} \cite[Lemma 3]{chan2025classificationclassplanarquadrinomials} Let $q$ be odd and $D(X)=\alpha X^2+\beta X+\alpha^q$ with $\alpha \in \gf_{q^2}$ and $\beta \in \gf_q$, not both zero. Define 
$\triangle(D):=\beta^2-4\alpha^{q+1}$. Then the following hold: 
\begin{itemize}
    \item[1)] $D(X)$ has a multiple root (which must be in $H$) if and only if $\triangle(D)=0$;
    \item[2)] $D(X)$ has two distinct roots in $H$ if and only if $\triangle(D)$ is a non-square in $\gf_q^*$;
    \item[3)] $D(X)$ has no roots in $H$ if and only if $\triangle(D)$ is a square in $\gf_q^*$. 
\end{itemize}
    
\end{lemma}

\begin{lemma}\cite{BRS67}\label{rootinfq}
	For a positive integer $n$, the quadratic equation $x^2+ax+b=0$, $a,b\in\gf_{2^n}$, $a\neq 0$ has solutions in $\gf_{2^n}$ if and only if $\mathrm{Tr}_{\gf_{2^n}/\gf_2}\left({b}/{a^2}\right)=0$. 
\end{lemma}

\begin{lemma}\cite[Lemma 2.4]{MR4626710}\label{rootinh2}
If $q$ is even and $A(X)=\alpha X^2+\beta X+\alpha^q$ with $\alpha \in \gf_{q^2}$ and $\beta \in \gf_{q}^*$, then the following are 
equivalent:

(1) $A(X)$ has at least one root in $H$; 

(2) $A(X)$ has two distinct roots in $H$; 

(3) $\mathrm{Tr}_{\gf_q/\gf_2} \left(\alpha^{q+1}/\beta^2\right)=1$. 
\end{lemma}

\section{Minimum distance of $C_{s}(q_0)$ and $\widetilde{C_{s}(q_0)}$}\label{sec3}

In this section we determine the minimum distance of $C_s(q_0)$ when $q_0$ is even and the minimum distance of $\widetilde{C_s(q_0)}$ which exists only when $q_0$ is odd. If $q_0$ is odd, it is trivial that the minimum distance of $C_s(q_0)$ is 2. 

\subsection{Case: $q_0$ is even}
When $q_0=2$, it was well-known that (see for example \cite{DHNX})
\begin{eqnarray*}
d(C_s(2))=\left\{\begin{array}{ll}
5,& \mbox{ if } s \mbox{ is even};\\
3,& \mbox{ if } s \mbox{ is odd}. 
\end{array}\right.
    \end{eqnarray*}
For $q_0\geq 4$ even, the result on $d\left(C_s(q_0)\right)$ is elementary but it is not easy for us to find a reference, so we compute the value $d(C_s(q_0))$ below.  
\begin{theorem} \label{x3:thm1} Let $q_0\geq 4$ be even. Then 	
	$$d(C_s(q_0)) = \left\{\begin{array}{ll}
		4, & \mbox{ if } s \mbox{ is even}; \\
		3, & \mbox{ if } s \mbox{ is odd}.
	\end{array} \right.$$
\end{theorem}
\begin{proof} First, recall that $C_s(q_0)$ is a $q_0$-ary $[q+1,q+1-2s]$ linear code and $d(C_s(q_0)) \le 4$ by \eqref{x:spbound}. Next, we establish a lower bound on the minimum distance. It is obvious that there is no codeword of weight 1 in $C_s(q_0)$. If there is a codeword of weight $2$ in $C_s(q_0)$, then there are $c_i,c_j\in\gf_{q_0}^*$, $0\leq i<j\leq q$, such that
\[c_i\xi^i+c_j\xi^j=0.\]
This implies $\xi^{j-i}=\frac{c_i}{c_j}\in\gf_{q_0}^* \cap H$. Since $\gf_{q_0}\cap H=\{1\}$, we have $\xi^{j-i}=1$, and consequently $j=i$, which is a contradiction. Hence, $d(C_s(q_0))\geq 3$ and thus we have $d(C_s(q_0)) \in \{3,4\}$. 

Assume now that $s$ is even. Suppose there exists a codeword of weight $3$ in $C_s(q_0)$, then there are $c_i,c_j,c_k\in\gf_{q_0}^*$ where $0\leq i<j<k\leq q$ such that
\begin{equation}\label{dijk1}
c_i\xi^i+c_j\xi^j+c_k\xi^k=0.
\end{equation}
Noting $ \xi^q = \xi^{-1} $, taking the $q$-th power on both sides of  Equation (\ref{dijk1}), we obtain
\begin{equation}\label{dijk2}
c_i\xi^{-i}+c_j\xi^{-j}+c_k\xi^{-k}=0.
\end{equation}
From (\ref{dijk1}) and (\ref{dijk2}), we derive
\[c^2_i=(c_i\xi^i)(c_i\xi^{-i})=\left(c_j\xi^j+c_k\xi^k\right)\left(c_j\xi^{-j}+c_k\xi^{-k}\right),\]
so we have 
\[\xi^{k-j}+\xi^{j-k}=\frac{c_i^2+c_j^2+c_k^2}{c_jc_k}.\]
Let $t=\xi^{k-j}$. The above equation becomes a quadratic equation in the variable $t$:
\begin{equation}\label{t}
t^2+\frac{c_i^2+c_j^2+c_k^2}{c_jc_k} \ t+1=0.
\end{equation}
There are two subcases. If $c_i^2+c_j^2+c_k^2=0$, then $t=1$, which implies $k=j$, a contradiction. If $c_i^2+c_j^2+c_k^2\neq 0$, by Lemma \ref{rootinfq}, Equation (\ref{t}) is solvable for $t \in \gf_{q}$ if and only if  
\[\mathrm{Tr}_{\gf_q/\gf_2}\left(\frac{(c_jc_k)^2}{(c_i^2+c_j^2+c_k^2)^2}\right)=0.\]
Since $s$ is even and $\frac{(c_jc_k)^2}{(c_i^2+c_j^2+c_k^2)^2} \in \gf_{q_0}$, the above condition always holds, so Equation (\ref{t}) is always solvable for $t\in\gf_q$. Thus, $t=\xi^{k-j}\in \gf_q\cap H=\{1\}$ and we have $k=j$, which is again a contradiction. Therefore, no codewords of weight $3$ exist in $C_s(q_0)$ and thus $d(C_s(q_0))=4$ when $s$ is even. 

Now we assume that $s$ is odd. Define 
$$H_{0}=\left\{x\in\gf_{q^2_0}:x^{q_0+1}=1\right\}.$$ 
Since $s$ is odd, $(q_0+1)|(q+1)$, $H_{0}$ is a subgroup of $H$. Let $\theta$ be a generator of $H_{0}$ and consider distinct elements $\theta^i,\theta^j \in H_{q_0 + 1} \setminus \{1\} $ with some $ 1 \leq i < j \leq q_0 $. It can be easily checked that
\[\frac{\theta^i(\theta^j+1)^2}{(\theta^i+\theta^j)^2},\quad \frac{\theta^j(\theta^i+1)^2}{(\theta^i+\theta^j)^2}\in\gf_{q_0}^*,\]
and 
\[\frac{\theta^i(\theta^j+1)^2}{(\theta^i+\theta^j)^2}\cdot \theta^i+\frac{\theta^j(\theta^i+1)^2}{(\theta^i+\theta^j)^2}\cdot \theta^j+1\cdot 1=0.\]
Since $1,\theta^i,\theta^j$ are distinct elements of $H_0$ and $H_{0}\subseteq H$, this gives a codeword of weight $3$ in $C_s(q_0)$. Hence, $d(C_s(q_0))=3$ when $s$ is odd. 
\end{proof}

\subsection{Case: $q_0$ is odd}

In the following, we investigate the minimum distance of $\widetilde{C_s(q_0)}$. For $q_0=3$, it was known from  \cite{D85,GS86} that $d\left(\widetilde{C_s(3)}\right) =5$ for any $s \ge 2$ (when $s=1$, the code $\widetilde{C_1(3)}$ is a trivial code of dimension 0, containing only the zero vector). For $q_0\geq 5$, We have the following theorem.


\begin{theorem} \label{x3:thm2} Let $q_0\geq 5$ be odd. Then 
		$$d\left(\widetilde{C_s(q_0)}\right) = \begin{cases}
		4, & ~\text{if}~s~\text{is~even} \\
		3, & ~\text{if}~s~\text{is~odd}.
	\end{cases} $$
\end{theorem}
\begin{proof} It is clear by construction that $d\left(\widetilde{C_s(q_0)}\right) \ge 3$. Now we assume that $s$ is even. We will show that there is no codeword in $\widetilde{C_s(q_0)}$ of weight $3$. Otherwise, by the definition of $\widetilde{C_s(q_0)}$, there exist $c_i,c_j,c_k\in\gf_{q_0}^*$, such that 
\begin{equation}\label{tildeeqn1}
	c_i\xi^i+c_j\xi^j+c_k\xi^k=0,
\end{equation}
where $0\leq i<j<k\leq \frac{q-1}{2}$.
Let $c_1=\frac{c_j}{c_i}$, $c_2=\frac{c_k}{c_i}$, $j'=j-i$ and $k'=k-i$. We have
\[1+c_1\xi^{j'}+c_2\xi^{k'}=0.\]
Taking the $q$-th power of the above equation, we have
\[1+c_1\xi^{-j'}+c_2\xi^{-k'}=0.\]
Hence,
\[c_2^2=c_2\xi^{k'}\cdot c_2\xi^{-k'}=(-c_1\xi^{j'}-1)(-c_1\xi^{-j'}-1),\]
i.e.,
\[(\xi^{j'})^2-\frac{c^2_2-c_1^2-1}{c_1}\xi^{j'}+1=0.\]
Consider the above equation as a quadratic equation in the variable $\xi^{j'}$. Then the discriminant is 
\[\Delta=\left(\frac{c^2_2-c_1^2-1}{c_1}\right)^2-4,\]
which is an element of $\gf_{q_0}$. Since $s$ is even, $\Delta$ is $0$ or an square element in $\gf_{q}$. Therefore, $\xi^{j'}\in\gf_{q}$. We conclude that $\xi^{j'}=1$ since $\gf_q\cap H_0=\{1\}$, which contradicts  $i<j$. Then $d(\widetilde{C_s(q_0)})\geq 4$ when $s$ is even. The equality holds because of \eqref{x:spbound}. 

In the following, we assume that $s$ is odd. We will prove the existence of a weight-3 codeword in  $\widetilde{C_s(q_0)}$. 

By Lemma \ref{exist}, there exist $c_1,c_2\in\gf_{q_0}\setminus\{0,\pm 1\}$ such that 
\[\delta=(c_1+c_2+1)(c_1+c_2-1)(c_1-c_2+1)(c_1-c_2-1)\]
is a nonsquare element in $\gf_{q_0}$. Since $s$ is odd, $\delta$ is a nonsquare element in $\gf_q$ as well. Let $\sqrt{\delta}$ be one of the square roots of $\delta$, which is in $\gf_{q^2}$. Define
\[\zeta_1=\frac{c^2_2-c^2_1-1+\sqrt{\delta}}{2c_1}, \quad \zeta_2=\frac{c^2_1-c^2_2-1-\sqrt{\delta}}{2c_2}.\]
One can verify that $\zeta_1$, $\zeta_2$ are two distinct elements in $H\setminus\{\pm 1\}$ and $\zeta_1\neq\zeta_2$. Then we can choose suitable $\epsilon_1,\epsilon_2\in\{\pm 1\}$ such that $\epsilon_1\zeta_1,\epsilon_2\zeta_2\in H_0\setminus\{1\}$ and they are distinct. Moreover, by a direct calculation, we have
\[(\epsilon_1c_1)\cdot(\epsilon_1\zeta_1)+(\epsilon_2c_2)\cdot (\epsilon_2\zeta_2)+1\cdot 1=0,\]
where $\epsilon_1c_1,\epsilon_2c_2\in\gf_{q_0}$, $\epsilon_1\zeta_1,\epsilon_2\zeta_2\in H_0\setminus\{1\}$ and they are distinct. Hence, $\widetilde{C_s(q_0)}$ has a codeword of weight $3$. We conclude that $d(\widetilde{C_s(q_0)})=3$ when $s$ is odd.  This completes the proof.

\end{proof} 

\section{Covering radius of $C_s(q_0)$ is at most 3}\label{sec4}

Denote by $\rho\left(C_s(q_0)\right)$ the covering radius of the code $C_s(q_0)$. In this section, we prove that $\rho\left(C_s(q_0)\right) \le 3$ for all $q_0$ and $s$. We first prove a technical lemma. 


\begin{lemma} \label{x:count}
For any $\alpha \in \gf_{q^2}$ and any $(c_1,c_2,c_3) \in \gf_{q_0}^3$, denote by $N_1$ the number of solutions $(x_1,x_2,x_3) \in H^3$ to the equation
\begin{equation}\label{x:main}
	\alpha=c_1x_1+c_2x_2+c_3x_3.
\end{equation}
Denote by $N_2$ the number of solutions $(x_1,x_2) \in H^2$ to the equation
\begin{equation}\label{x:thm2-2}
  c_2(\overline{\alpha}-c_1\overline{x}_1)x^2_2-\left((\alpha-c_1x_1)(\overline{\alpha}-c_1\overline{x}_1)+c^2_2-c^2_3\right)x_2+c_2(\alpha-c_1x_1)=0.
 \end{equation}
If $c_3 \ne 0$, then $N_1=N_2$. 
\end{lemma}
\begin{proof}
Suppose $(x_1,x_2,x_3) \in H^3$ is a solution to Eq \eqref{x:main}, then we have 
\[c_3x_3=\alpha-c_1x_1-c_2x_2.\]
Since $c_i\in \gf_{q_0} \subseteq \gf_q$ for each $i$, we have 
 \[c_3\overline{x}_3=\overline{\alpha}-c_1\overline{x}_1-c_2\overline{x}_2.\]
Here for $x \in \gf_{q^2}$ we denote $\overline{x}:=x^q$. Since $x_3 \in H$, $\overline{x}_3x_3=1$, so we have 
 \begin{eqnarray} \label{x:c3} c_3^2=(\alpha-c_1x_1-c_2x_2)(\overline{\alpha}-c_1\overline{x}_1-c_2\overline{x}_2).\end{eqnarray}
From this we obtain 
 \begin{equation}\label{x:thm2-1}
 (\alpha-c_1x_1)(\overline{\alpha}-c_1\overline{x}_1)-c_2x_2(\overline{\alpha}-c_1\overline{x}_1)-c_2\overline{x}_2(\alpha-c_1x_1)+c_2^2=c^2_3.
 \end{equation}
Multiplying $x_2$ on both sides of Eq \eqref{x:thm2-1} and noting $x_2\overline{x}_2=1$, we obtain Eq \eqref{x:thm2-2}. 

On the other hand, suppose $(x_1,x_2) \in H^2$ is a solution to Eq \eqref{x:thm2-2}, then $(x_1,x_2)$ also satisfies Eq \eqref{x:thm2-1} and hence Eq \eqref{x:c3}, so $\alpha-c_1x_1-c_2x_2 \ne 0$ as $c_3 \ne 0$. It is easy to see that by defining  
\[x_3:=\frac{\alpha-c_1x_1-c_2x_2}{c_3},\]
we have $x_3 \in H$ and this triple $(x_1,x_2,x_3)\in H^3$ is a solution to Eq \eqref{x:main}. This completes the proof of Lemma \ref{x:count}. 
\end{proof}

We now proceed to prove that $\rho\left(C_s(q_0)\right) \le 3$. While the proofs for even and odd $q_0$ are similar in nature, the technical treatments are different and non-trivial. For the sake of clarity, we write these two cases as two separate theorems. It is worth pointing out that the result for odd $q_0$ was already obtained (see \cite[Theorem II.1]{SHO23} and \cite[Corollary 3.14]{SLHO25}). Here for odd $q_0$, we provide a new proof. 

\begin{theorem} \label{x2:thmodd} Let $ q_0 \geq 3 $ be odd and $ s \geq 1 $. Then the covering radius of $C_{s}(q_0)$ is at most $3$. 
\end{theorem}

\begin{proof}
The case that $s=1$ is trivial as was noted in \cite{SHO23}. From now on let us assume that $s \geq 2$. To prove that $\rho(C_{s}(q_0))\leq 3$, by Lemma \ref{x:crlemma}, it suffices to show that for any $\alpha\in\gf_{q^2}$, there exist $c_1,c_2,c_3\in \gf_{q_0}$ and $x_1,x_2,x_3\in H$, such that 
\begin{equation}\label{sgeq2}
\alpha=c_1x_1+c_2x_2+c_3x_3.
\end{equation}
If $\alpha\in \gf_{q_0}\cdot H=\left\{c x:c\in\gf_{q_0}, x\in H\right\}$, then taking  $c_1\in \gf_{q_0}$ and $x_1\in H$ such that $\alpha=c_1x_1$ and taking $c_2=c_3=0$, we see that (\ref{sgeq2}) holds for these $c_i$ and $x_i$. Now suppose $\alpha\in \gf_{q^2}\setminus\gf_{q_0}\cdot H$. For any fixed $(c_1,c_2,c_3)\in \gf_{q_0}^{*3}$, let us define 
 \[N(c_1,c_2,c_3):=\#\left\{(x_1,x_2,x_3)\in H^3: \alpha=c_1x_1+c_2x_2+c_3x_3\right\}.\]
We will prove that $N(c_1,c_2,c_3)>0$ for any fixed $(c_1,c_2,c_3)\in \gf_{q_0}^{*3}$ and any $\alpha\in \gf_{q^2}\setminus\gf_{q_0}\cdot H$. This would imply that $\rho(C_{s}(q_0))\leq 3$. 
 
By Lemma \ref{x:count}, $N(c_1,c_2,c_3)$ equals the number of solutions $(x_1,x_2)\in H^2$ to the equation 
 \begin{equation}\label{thm2-2}
  c_2(\overline{\alpha}-c_1\overline{x}_1)x^2_2-\left((\alpha-c_1x_1)(\overline{\alpha}-c_1\overline{x}_1)+c^2_2-c^2_3\right)x_2+c_2(\alpha-c_1x_1)=0.
 \end{equation}
Since $\alpha \in \gf_{q^2} \setminus \gf_{q_0} \cdot H$, $\alpha-c_1x_1\neq 0$, hence for any fixed $x_1 \in H$, Eq (\ref{thm2-2}) is a quadratic equation in $x_2$, whose discriminant is given by 
\begin{eqnarray*} \Delta(x_1)&=&\left((\alpha-c_1x_1)(\overline{\alpha}-c_1\overline{x}_1)+c_2^2-c_3^2\right)^2-4c^2_2(\alpha-c_1x_1)(\overline{\alpha}-c_1\overline{x}_1). 
\end{eqnarray*}
It is also easy to check that 
\begin{eqnarray} \Delta(x_1)
&=&\left((\alpha-c_1x_1)(\overline{\alpha}-c_1\overline{x}_1)-(c_2+c_3)^2\right)\left((\alpha-c_1x_1)(\overline{\alpha}-c_1\overline{x}_1)-(c_2-c_3)^2\right). \label{x:deltaodd} 
\end{eqnarray}
According to Lemma \ref{x:chan-X}, for any given $x_1 \in H$, the number of roots $x_2 \in H$ satisfying Eq \eqref{thm2-2} is exactly $1-\chi(\Delta(x_1))$, where $\chi$ denotes the quadratic character of $\gf_{q}$, so we derive that 
\begin{align*}
	 N(c_1,c_2,c_2)&=\sum_{x_1\in H}\left(1-\chi(\Delta(x_1))\right)
	 =q+1-\sum_{x_1\in H}\chi(\Delta(x_1)).
\end{align*}

Now we fix an element $a \in \mathbb{F}_{q^2} \setminus \mathbb{F}_q$ and consider the bijection $\phi: \mathbb{P}^1(\gf_q) \to H$ given by $\phi: y \mapsto \frac{y+a}{y+\overline{a}}$. Since $\phi(\infty)=1 \in H$, we have 
\begin{align*}
	 N(c_1,c_2,c_2)
	 &=q+1-\chi\left(\triangle(1)\right)-\sum_{x\in\gf_q}\chi\left(\Delta\left(\frac{x+a}{x+\overline{a}}\right)\right).
\end{align*}
For the function $\Delta(\frac{x+a}{x+\overline{a}})$, from \eqref{x:deltaodd} we can obtain 
 \begin{eqnarray*}
 \Delta\left(\frac{x+a}{x+\overline{a}}\right)&=&\left(\left(\alpha-c_1\frac{x+a}{x+\overline{a}}\right)\left(\overline{\alpha}-c_1\frac{x+\overline{a}}{x+a}\right)-(c_2+c_3)^2\right) \times \\
 & &\left(\left(\alpha-c_1\frac{x+a}{x+\overline{a}}\right)\left(\overline{\alpha}-c_1\frac{x+\overline{a}}{x+a}\right)-(c_2-c_3)^2\right)\\
 &=&\left(\frac{(\alpha-c_1)(\overline{\alpha}-c_1)}{(x+a)(x+\overline{a})}\right)^2f_1(x)f_2(x),
 \end{eqnarray*}
where 
\begin{eqnarray*} 
f_1(x)&=&\left(x+\frac{\alpha \overline{a}-c_1a}{\alpha-c_1}\right)\left(x+\frac{\overline{\alpha} {a}-c_1\overline{a}}{\overline{\alpha}-c_1}\right)-\frac{(c_2+c_3)^2}{(\alpha-c_1)(\overline{\alpha}-c_1)}, \\
f_2(x)&=&\left(x+\frac{\alpha \overline{a}-c_1a}{\alpha-c_1}\right)\left(x+\frac{\overline{\alpha} {a}-c_1\overline{a}}{\overline{\alpha}-c_1}\right)-\frac{(c_2-c_3)^2}{(\alpha-c_1)(\overline{\alpha}-c_1)}. 
\end{eqnarray*}
Since $a \in \gf_{q^2} \setminus \gf_q$, $c_1 \in \gf_{q_0}^*$ and $\alpha \in \gf_{q^2} \setminus \gf_{q_0} \cdot H$, we have
\[\frac{(\alpha-c_1)(\overline{\alpha}-c_1)}{(x+a)(x+\overline{a})} \in \gf_q^*, \quad \forall x \in \gf_q,\]
so 
\[\chi\left(\Delta\left(\frac{x+a}{x+\overline{a}}\right)\right)=\chi\left(f_1(x)f_2(x)\right), \quad \forall x \in \gf_{q},\]
and thus
\[N(c_1,c_2,c_3)= q+1-\chi\left(\triangle(1)\right)-\sum_{x \in \gf_q} \chi\left(f_1(x)f_2(x)\right).\]
As $c_2,c_3 \in \gf_{q_0}^*$, $(c_2+c_3)^2 \ne (c_2-c_3)^2$, it is easy to see that $\gcd(f_1(X),f_2(X))=1$, and the quadratic polynomials $f_1(X)$ and $f_2(X)$ can not be squares of linear polynomials simultaneously, hence by Lemma \ref{weil}, we have 
\[\left|\sum_{x\in\gf_q}\chi\left(f_1(x)f_2(x)\right)\right|\leq 3\sqrt{q}.\]
This gives 
\[N(c_1,c_2,c_3)\ge q-\left|\sum_{x\in\gf_q}\chi\left(f_1(x)f_2(x)\right)\right|\geq q-3\sqrt{q}>0,\]
as long as $q>9$. Hence we obtain $\rho\left(C_s(q_0)\right)\leq 3$ if $q>9$.
 
For $q \le 9$ and odd, since $s \ge 2$, the only possibility is $q_0=3,s=2$ and $q=9$. It can be easily checked by {\sc Magma} that $\rho\left(C_2(3)\right)=3$. This completes the proof of Theorem \ref{x2:thmodd}. 
\end{proof}

\begin{theorem} \label{x2:thm} Let $ q_0 \geq 2 $ be even and $ s \geq 1 $. Then the covering radius of $C_{s}(q_0)$ is at most $3$. 
\end{theorem}
\begin{proof}

The case $s=1$ is trivial, because in this case $q=q_0$, by the standard polar representation, every $x\in\gf_{q_0^2}^*$ can be written uniquely as $x=yz$ with $y\in\gf_{q_0}^*$ and $z\in H$. Therefore in this case $\rho(C_{1}(q_0))=1$.

If $q_0=2$ and $s=2$, it is easy to check by {\sc Magma} that $\rho(C_s(2))=2$. So from now on let us assume either  ``$q_0=2$ and $s \ge 3$'' or ``$q_0 \ge 4$ and $s \ge 2$'', so we always have $q \ge 8$. 

By Lemma \ref{x:crlemma}, to show that $\rho(C_s(q_0)) \leq 3$, it suffices to prove that for every $\alpha \in \mathbb{F}_{q^2}$, there exist $c_1, c_2, c_3 \in \mathbb{F}_{q_0}$ and $x_1, x_2, x_3 \in H$ such that
\begin{equation}\label{main}
	\alpha=c_1x_1+c_2x_2+c_3x_3.
\end{equation}
We may write $\alpha$ as $\alpha = \alpha_0 \alpha_H$ with $\alpha_0 \in \mathbb{F}_q$ and $\alpha_H \in H$. Suppose~\eqref{main} holds for $\alpha$. Then
\[
\alpha_0 = c_1 \frac{x_1}{\alpha_H} + c_2 \frac{x_2}{\alpha_H} + c_3 \frac{x_3}{\alpha_H},
\]
where $\frac{x_i}{\alpha_H} \in H$ for each $i$ since $H$ is the multiplicative subgroup of order $q+1$ in $\gf_{q^2}^*$. Therefore, without loss of generality, we may assume $\alpha \in \mathbb{F}_q$. If $\alpha \in \mathbb{F}_{q_0}$, we can set $c_1 = \alpha$, $x_1 = 1$, and $c_2 = c_3 = 0$, so~\eqref{main} holds trivially. Hence, from now on let us assume that $\alpha \in \mathbb{F}_q \setminus \mathbb{F}_{q_0}$.

For fixed $(c_1, c_2, c_3) \in \mathbb{F}_{q_0}^{*3}$, let us define
\[
N(c_1, c_2, c_3) = \# \left\{ (x_1, x_2, x_3) \in H^3 : \alpha = c_1 x_1 + c_2 x_2 + c_3 x_3 \right\}.
\]
Our goal is to show that $N(c_1, c_2, c_3) >0$ for any $(c_1, c_2, c_3) \in (\mathbb{F}_{q_0}^*)^3$ and any $\alpha \in \gf_q \setminus \gf_{q_0}$. This would imply that $\rho\left(C_s(q_0)\right) \leq 3$. 

By Lemma \ref{x:count} again, $N(c_1,c_2,c_3)$ equals the number of solutions $(x_1,x_2) \in H^2$ to the equation 
\begin{equation}\label{thm2-3}
	c_2\left(\alpha+c_1\overline{x}_1\right)x_2^2+\left((\alpha+c_1x_1)(\alpha+c_1\overline{x}_1)+c_2^2+c_3^2\right)x_2+c_2\left(\alpha+c_1x_1\right)=0.
\end{equation}
For any fixed $x_1 \in H$, to analyze the existence of solutions $x_2 \in H$ of Eq (\ref{thm2-3}), we denote
\begin{eqnarray*}
    \triangle(x_1):=(\alpha+c_1x_1)(\alpha+c_1\overline{x}_1)+c_2^2+c_3^2. \label{x:tri-even}
\end{eqnarray*}
We consider two cases.

{\bf Case 1:  $\triangle(x_1)=0$.} Then~\eqref{thm2-3} reduces to
\[(\alpha+c_1\overline{x}_1)c_2x^2_2+(\alpha+c_1x_1)c_2=0,\]
that is,
\[x_2^2=\frac{\alpha+c_1x_1}{\alpha+c_1\overline{x}_1}.\]
Since $q$ is even, this gives a unique solution $x_2\in H$ to Eq \eqref{thm2-3}. 

{\bf Case 2: $\triangle(x_1) \ne 0$.} Then by Lemma~\ref{rootinh2}, for each fixed $x_1 \in H$, Eq (\ref{thm2-3}) has one and hence two distinct solutions $x_2 \in H$ if and only if
\begin{equation}\label{thm2-4}
	\mathrm{Tr}_{\gf_q/\gf_2}\left(\frac{c_2^2(\alpha+c_1\overline{x}_1)(\alpha+c_1x_1)}{\left((\alpha+c_1x_1)(\alpha+c_1\overline{x}_1)+c_2^2+c_3^2\right)^2}\right)=1.
\end{equation}
Note that the left hand side can be rewritten as 
\begin{align*}
	&\mathrm{Tr}_{\gf_q/\gf_2}\left(\frac{c_2^2\left((\alpha+c_1\overline{x}_1)(\alpha+c_1x_1)+c_2^2+c_3^2\right)+c^2_2\left(c_2^2+c_3^2\right)}{\left((\alpha+c_1x_1)(\alpha+c_1\overline{x}_1)+c_2^2+c_3^2\right)^2}\right)\\	=&\mathrm{Tr}_{\gf_q/\gf_2}\left(\frac{c_2^2}{(\alpha+c_1x_1)(\alpha+c_1\overline{x}_1)+c_2^2+c_3^2}\right)+\mathrm{Tr}_{\gf_q/\gf_2}\left(\frac{c^2_2\left(c_2^2+c_3^2\right)}{\left((\alpha+c_1x_1)(\alpha+c_1\overline{x}_1)+c_2^2+c_3^2\right)^2}\right)\\
	=&\mathrm{Tr}_{\gf_q/\gf_2}\left(\frac{c_2^2}{(\alpha+c_1x_1)(\alpha+c_1\overline{x}_1)+c_2^2+c_3^2}\right)+\mathrm{Tr}_{\gf_q/\gf_2}\left(\frac{c_2(c_2+c_3)}{(\alpha+c_1x_1)(\alpha+c_1\overline{x}_1)+c_2^2+c_3^2}\right)\\	=&\mathrm{Tr}_{\gf_q/\gf_2}\left(\frac{c_2c_3}{(\alpha+c_1x_1)(\alpha+c_1\overline{x}_1)+c_2^2+c_3^2}\right). 
\end{align*}
So the following condition is equivalent to (\ref{thm2-4}): 
\begin{equation*}\label{thm2-5}
\mathrm{Tr}_{\gf_q/\gf_2}\left(\frac{c_2c_3}{(\alpha+c_1x_1)(\alpha+c_1\overline{x}_1)+c_2^2+c_3^2}\right)=1.
\end{equation*}
Hence for the given $x_1 \in H$, the number of solutions $x_2  \in H$ to Eq (\ref{thm2-3}) is given by 
\[1-\psi\left(\frac{c_2c_3}{(\alpha+c_1x_1)(\alpha+c_1\overline{x}_1)+c_2^2+c_3^2}\right),\]
where $\psi: \beta \mapsto (-1)^{\mathrm{Tr}_{\gf_q/\gf_2}(\beta)}$ is the standard additive character defined on $\gf_{q}$. 
Summarizing {\bf Cases 1} and {\bf 2} above, we conclude that 
\begin{eqnarray*} \label{x3:N}
N(c_1,c_2,c_3)&=&\sum_{\substack{x_1 \in H\\
\triangle(x_1)=0}} 1+\sum_{\substack{x_1 \in H\\
\triangle(x_1) \ne 0}} \left(1-\psi\left(\frac{c_2c_3}{(\alpha+c_1x_1)(\alpha+c_1\overline{x}_1)+c_2^2+c_3^2}\right)\right)\\
&=&q+1-\sum_{\substack{x_1 \in H\\
\triangle(x_1) \ne 0}}\psi\left(\frac{c_2c_3}{(\alpha+c_1x_1)(\alpha+c_1\overline{x}_1)+c_2^2+c_3^2}\right).
\end{eqnarray*}
Let us define 
\[
T := \sum_{\substack{x_1 \in H\\
\triangle(x_1) \ne 0}} \psi\left( \frac{c_2 c_3}{ (\alpha + c_1 x_1)(\alpha + c_1 \overline{x}_1) + c_2^2 + c_3^2 } \right). 
\]
We fix an element $a \in \mathbb{F}_{q^2} \setminus \mathbb{F}_q$ and consider the bijection $\phi: \mathbb{P}^1(\gf_q) \to H$ given by $\phi: y \mapsto \frac{y+a}{y+\overline{a}}$. Since $\phi(\infty)=1 \in H$ and 
\[\triangle(1)=\left(\alpha+c_1+c_2+c_3\right)^2 \ne 0,\]
using this bijection $\phi$ we can rewrite $T$ as 
\[T=\psi\left( \frac{c_2 c_3}{ \left(\alpha + c_1 + c_2 + c_3\right)^2 } \right)+\sum_{\substack{y \in \gf_{q}\\
\triangle\left(\phi(y)\right) \ne 0}} \psi\left( \frac{c_2 c_3}{ \left(\alpha + c_1 \phi(y)\right)\left(\alpha + c_1 \phi(y)^{-1}\right) + c_2^2 + c_3^2 } \right).\]
The right hand side can be rewritten as 
\[
T =\psi\left(\frac{c_2c_3}{(\alpha+c_1+c_2+c_3)^2}\right)+\sum_{y \in \mathbb{F}_q \setminus S} \psi\left( f(y) \right),
\]
where $f(y)$ is a rational function in $\mathbb{F}_q(y)$ given by 
\[f(y)=\frac{c_2c_3}{\left(\alpha+c_1+c_2+c_3\right)^2} \left(1+\frac{A}{y^2+\left(a+\overline{a}\right)y+a\overline{a}+A}\right),\]
$S$ is the set of poles of $f(y)$ in $\gf_q$ and  
\[A=\frac{\alpha c_1\left(a+\overline{a}\right)^2}{\left(\alpha+c_1+c_2+c_3\right)}.\]
Now we apply Lemma~\ref{rational}. From the explicit form of $ f(y) $, it is easy to see that $M=0$ and $L=4$, so we deduce
\[
|T| \leq 1+ \left(1 + 2\sqrt{q}\right)=2 \left(1+\sqrt{q}\right).
\]
Armed with this estimate, we obtain 
\[N(c_1,c_2,c_3) \ge q+1-|T| \ge q+1-2 \left(1+\sqrt{q}\right)=q-2\sqrt{q}-1.
\]
Since we have assumed that $q \geq 8$, it is easy to check that $q-2\sqrt{q}-1 > 0$, so $N(c_1,c_2,c_3)>0$ as desired for any $(c_1,c_2,c_3) \in \gf_{q_0}^{*3}$ and any $\alpha \in \gf_q \setminus \gf_{q_0}$. This completes the proof of Theorem \ref{x2:thm}. 
\end{proof}
Since $\rho\left(C_s(q_0)\right) >1$ if $s \ge 2$, we summarize Theorems \ref{x2:thmodd} and \ref{x2:thm} in the following corollary.
\begin{corollary} \label{x3:cor} For any $q_0$, if $s\geq 2$, then the covering radius of $ C_s(q_0) $ satisfies $ \rho\left(C_s(q_0)\right) \in \{2, 3\} $.
\end{corollary}

\section{Determining covering radius of $C_{s}(q_0)$: $q_0$ odd}\label{sec5}

Let $ q_0 \geq 3 $ be odd and $ s \geq 2 $. It is known by Corollary \ref{x3:cor} that $\rho\left(C_s(q_0)\right) \in \{2,3\}$. In this section we study explicit conditions to distinguish these two cases. While this case was treated in \cite{SHO23,SLHO25}, we will provide new proofs and stronger results. 

We know by the proof of Theorem \ref{x2:thmodd} that $ \rho(C_s(q_0)) = 2 $ if and only if for any $ y \in \gf_{q^2} \setminus \gf_{q_0} \cdot H $, there exist $(c_1, c_2) \in \mathbb{F}_{q_0}^{*2}$ and $ (x_1, x_2) \in H^2 $ such that 
\begin{eqnarray} \label{x5:r2odd}
y = c_1 x_1 + c_2 x_2.
\end{eqnarray}
Taking the $ q $-th power on both sides of the above equation, using $\overline{x}=x^q$ for any $x \in \gf_{q^2}$ we obtain
$$
\overline{y} = c_1 \overline{x}_1 + c_2 \overline{x}_2.
$$
Hence, $ x_1 $ satisfies the equation
$$
\left(y-c_1 x_1\right)(\overline{y}-c_1 \overline{x}_1) =c_2^2x_2\overline{x}_2= c_2^2,
$$
which simplifies to
\begin{eqnarray} \label{x4:eq1odd}
\overline{y}x_1^2 - \left(\frac{y \overline{y} + c_1^2 - c_2^2}{c_1}\right) \, \, x_1 + y = 0.
\end{eqnarray}
For the given $y \in \gf_{q^2} \setminus \gf_{q_0}\cdot H$ and $c_1,c_2 \in \gf_{q_0}^*$, if the above quadratic equation has a solution $x_1 \in H$, then by taking $x_2=\frac{y-c_1x_1}{c_2}$, it is easy to see that $x_2 \in H$ and this pair $(x_1,x_2)\in H^2$ satisfies Eq (\ref{x5:r2odd}). So to study Eq (\ref{x5:r2odd}), we just need to study Eq (\ref{x4:eq1odd}). 

The discriminant of the quadratic equation \eqref{x4:eq1odd} is 
\begin{eqnarray*}  \triangle:=\left(\frac{y \overline{y} + c_1^2 - c_2^2}{c_1}\right)^2-4y \overline{y}, \end{eqnarray*}
and according to Lemma \ref{x:chan-X}, the number of solutions $x_1 \in H$ to Eq \eqref{x4:eq1odd} is exactly $1-\chi(\triangle)$, here $\chi$ is the quadratic character of $\gf_{q}$. So $\rho(C_s(q_0)=2$ if and only if for any $y \in \gf_{q^2} \setminus \gf_{q_0} \cdot H$, there exist $(c_1,c_2) \in \gf_{q_0}^{*2}$ such that 
\[1-\chi(\triangle)>0,\]
and thus $\rho(C_s(q_0)=3$ if and only if there exists $y \in \gf_{q^2} \setminus \gf_{q_0} \cdot H$ such that for any $(c_1,c_2) \in \gf_{q_0}^{*2}$ we have  
\[1-\chi(\triangle)=0 \Longleftrightarrow \triangle \in \square_q \subseteq \gf_{q}^*.\]

We summarize the above in the following lemma.

\begin{lemma}\label{rho=3-odd}
	Let $ q_0 $ be odd and $ s \geq 2 $. Then $ \rho(C_s(q_0)) = 3 $ if and only if there exists $y \in \gf_{q^2} \setminus \mathbb{F}_{q_0} \cdot H $ such that \begin{equation}\label{thm2-6-odd}
	\left(\frac{y \overline{y} + c_1^2 - c_2^2}{c_1}\right)^2-4y \overline{y} \in \square_{q}, \quad \forall c_1,c_2 \in \gf_{q_0}^*.
	\end{equation}
\end{lemma}

Let us simplify Eq (\ref{thm2-6-odd}) further in Lemma \ref{rho=3-odd}. It is easy to see that 
\begin{eqnarray*} \left(y \overline{y} + c_1^2 - c_2^2\right)^2-4c_1^2y \overline{y}
=\left(y \overline{y}-(c_1+c_2)^2\right)\left(y \overline{y}-(c_1-c_2)^2\right).\end{eqnarray*}
Setting \[\alpha=c_1+c_2, \quad \beta=c_1-c_2,\]
then 
\[c_1=\frac{1}{2}\left(\alpha+\beta\right), \quad c_2=\frac{1}{2}\left(\alpha-\beta\right).\]
It is easy to see that 
\[c_1,c_2 \in \gf_{q_0}^* \Longleftrightarrow \alpha \ne \pm \beta \Longleftrightarrow \alpha^2 \ne \beta^2, \alpha,\beta \in \gf_{q_0}.\] 
So Eq \eqref{thm2-6-odd} is equivalent to 
\[\left(y\overline{y}-\alpha^2\right)\left(y\overline{y}-\beta^2\right) \in \square_q, \quad \forall \alpha,\beta \in \gf_{q_0} \mbox{ with } \alpha^2 \ne \beta^2.\]
Considering the case that $\alpha=0$, the above can be further simplified as 
\[y\overline{y}\left(y\overline{y}-\beta\right) \in \square_q, \quad \forall \beta \in \square_{q_0}.\]
Finally, $y \in \gf_{q^2} \setminus \gf_{q_0} \cdot H$ if and only if $y\overline{y} \in \gf_{q}^* \setminus \square_{q_0}$. 
Summarizing the above results, we obtain the following simplified characterization.

\begin{theorem}
	\label{r3s-odd}
	Let $q_0$ be odd and $s \geq 2$. Then $\rho(C_s(q_0)) = 3$ if and only if there exists $x \in \gf_{q}^* \setminus \square_{q_0}$ such that 
\begin{eqnarray*} \label{x6:r3-odd}
x \left(x-\beta\right) \in \square_q, \quad \forall \beta \in \square_{q_0}.
\end{eqnarray*}
Otherwise $\rho(C_s(q_0))=2$.     
\end{theorem}

As consequences of Theorem~\ref{r3s-odd}, we collect below several results about $\rho(C_s(q_0))$ for some special $q_0$ and $s$.

\begin{corollary} \label{x4:cor2-odd} Let $q_0$ be odd. If $s \ge 3$ is odd and 
\begin{eqnarray}\label{x4:sodd}
s<1+\frac{1}{2}\left(\sqrt{q_0}-\frac{1}{\sqrt{q_0}}\right),\end{eqnarray}
then $\rho(C_s(q_0))=2$. 
\end{corollary}

\begin{proof}
Suppose $\rho(C_s(q_0))=3$, by Theorem \ref{r3s-odd}, there is an element $x_0 \in \gf_{q}^*\setminus \square_{q_0}$ such that 
\[x_0 \left(x_0-\beta\right) \in \square_q, \quad \forall \beta \in \square_{q_0}.\]
Let $t$ be the least positive integer such that $x_0 \in \gf_{q_0^t}$. Since $q=q_0^s$, obviously $t|s$. For any fixed $\beta \in \square_{q_0}$, since $x_0\left(x_0-\beta\right) \in \square_q$, we have 
\[\left(x_0\left(x_0-\beta\right)\right)^{q_0^i}=x_0^{q_0^i}\left(x_0^{q_0^i}-\beta\right) \in \square_q, \quad \forall 0 \le i \le t-1.  \]
Multiplying these terms  together for $0 \le i \le t-1$, also noting that $s$ is odd, we find that 
\begin{eqnarray} \label{x:normv} \mathrm{N}_{\gf_{q_0^t}/\gf_{q_0}} \left(x_0\left(x_0-\beta\right)\right)=\mathrm{N}_{\gf_{q_0^t}/\gf_{q_0}} (x_0) \mathrm{N}_{\gf_{q_0^t}/\gf_{q_0}} \left(x_0-\beta\right) \in \square_{q} \cap \gf_{q_0}=\square_{q_0}, \quad \forall \beta \in \square_{q_0}.\end{eqnarray}
Here $\mathrm{N}_{\gf_{q_0^t}/\gf_{q_0}}: \gf_{q_0^t} \to \gf_{q_0}$ is the Norm function given by $\mathrm{N}_{\gf_{q_0^t}/\gf_{q_0}}(x)=x \cdot x^{q_0} \cdots x^{q_0^{t-1}}$ for any $x \in \gf_{q_0^t}$. 

Now we evaluate this norm function explicitly. Let $f(X)$ be the monic irreducible polynomial of $x_0$ over $\gf_{q_0}$. Then $f(X) \in \gf_{q_0}[X], \deg f=t$ and  
\begin{eqnarray*} \label{x:fx} f(X)=\prod_{i=0}^{t-1} \left(X-x_0^{q_0^i}\right),\end{eqnarray*}
so we find  
\[f(0)=\left(-x_0\right) \cdot \left(-x_0^{q_0}\right) \cdots \left(-x_0^{q_0^{t-1}}\right)= (-1)^t \mathrm{N}_{\gf_{q_0^t}/\gf_{q_0}}(x_0). \]
Since $t|s$ and $s$ is odd, we have
\begin{eqnarray} \label{x:norm}
    \mathrm{N}_{\gf_{q_0^t}/\gf_{q_0}}(x_0)=-f(0).
\end{eqnarray}
Next, denote 
\[g_{\beta}(X):=f(X+\beta), \quad \forall \beta \in \square_{q_0}.\]
It is easy to see that $g_{\beta}(X)$ is the monic irreducible polynomial of $x_0-\beta$ over $\gf_{q_0}$. From Eq \eqref{x:norm} we obtain
\begin{eqnarray} \label{x:norm2} \mathrm{N}_{\gf_{q_0^t}/\gf_{q_0}}\left(x_0-\beta\right)=-g_{\beta}(0)=-f(\beta). \end{eqnarray}
Plugging \eqref{x:norm} and \eqref{x:norm2} into Eqs \eqref{x:normv} we find that 
\begin{eqnarray} \label{x:fs}
    f(0)f(\beta) \in \square_{q_0}, \qquad \forall \beta \in \square_{q_0}, 
\end{eqnarray}
where $f(X)\in \gf_{q_0}[X]$ is the monic irreducible polynomial of $x_0$ over $\gf_{q_0}$. 

Let $\chi$ denote the quadratic character defined on $\gf_{q_0}$. Eqs \eqref{x:fs} imply that 
\[\chi\left(f(x^2)\right)=\chi(f(0)) \in \left\{\pm 1 \right\}, \quad \forall x \in \gf_{q_0}. \]
So we have
\begin{eqnarray} \label{x:A1} A:=\left|\sum_{x \in \gf_{q_0}}\chi\left(f(x^2)\right)\right|=q_0. \end{eqnarray}
On the other hand, denote $g(X):=f(X^2)$. Let $\lambda \in \gf_{q^2}^*$ be such that $\lambda^2=x_0$, then $\left(\pm \lambda^{q_0^i}\right)^2=\lambda^{2q_0^i}=x_0^{q_0^i}$ for any $0 \le i \le t-1$, from \eqref{x:fx} we can obtain
\[g(X)=\prod_{i=0}^{t-1} \left(X^2-x_0^{q_0^i}\right)=\prod_{i=0}^{t-1} \left(X-\lambda^{q_0^i}\right)\left(X+\lambda^{q_0^i}\right).\]
It is easy to see that all the roots $\pm \lambda^{q_0^i}, 0 \le i \le t-1$ are distinct: suppose $\epsilon_1 \lambda^{q_0^i}=\epsilon_2 \lambda^{q_0^j}$ for some $\epsilon_1,\epsilon_2 \in \{\pm 1\}$ and $0 \le i,j \le t-1$. Taking squares on both sides, we obtain $x_0^{q_0^i}=x_0^{q_0^j}$, thus $i=j$ from which we also obtain $\epsilon_1=\epsilon_2$. So $g(X)$ is not a square of polynomials. By Lemma \ref{weil} we obtain 
\begin{eqnarray} \label{x:A2} A=\left|\sum_{x \in \gf_{q_0}}\chi\left(f(x^2)\right)\right|=\left|\sum_{x \in \gf_{q_0}}\chi\left(g(x)\right)\right|\le 1+(2t-2)\sqrt{q_0}.\end{eqnarray}
Combining \eqref{x:A1} and \eqref{x:A2}, we obtain
\[A=q_0 \le 1+(2t-2)\sqrt{q_0} \le 1+(2s-2)\sqrt{q_0}. \]
That is,
\[2(s-1) \ge \sqrt{q_0}-\frac{1}{\sqrt{q_0}}\Longrightarrow s \ge 1+\frac{1}{2} \left(\sqrt{q_0}-\frac{1}{\sqrt{q_0}}\right).\]
This contradicts the assumption that $s < 1+\frac{1}{2} \left(\sqrt{q_0}-\frac{1}{\sqrt{q_0}}\right)$. This completes the proof of Corollary \ref{x4:cor2-odd}. 
\end{proof}
For the sake of completeness, we also include the following results from \cite{SHO23}.
\begin{corollary} \label{x4:cor3-odd}
	Let $q_0$ be odd and $s \ge 2$. 
    \begin{itemize}
        \item[0).] If $s \ge 2$, then $\rho(C_s(3))=3$;
        \item[1).] If $s\geq 2$ is even, then $\rho(C_s(q_0))=3$;
        \item[2).] If $s$ is odd, and there is $s'|s$ such that $\rho(C_{s'}(q_0))=3$, then $\rho(C_s(q_0))=3$; 
    \end{itemize}
    \end{corollary}
\begin{proof}
Items 0) 1) and 2) were all proved in \cite{SHO23}. In particular they correspond to \cite[Corollary III.1, Corollary III.3 and Proposition V.1]{SHO23} respectively. We omit the details. 
\end{proof}

Finally, similar to \cite[Theorem IV.4]{SHO23} and \cite[Theorem 4.8]{SLHO25}, for any $q_0 \ge 3$ odd and $s$ odd, we prove that if $s$ is sufficiently large, then $\rho\left(C_s(q_0)\right)=3$. 

\begin{theorem} \label{x4:thm4-odd} Let $q_0 \ge 3$ be odd and $s \ge 3$ be odd. Denote $m:=\frac{q_0-1}{2}$. If 
\begin{eqnarray} \label{x:bounds-odd}
    q_0^s-q_0^{s/2}\left((m-3)2^{m-1}+2\right)>3 \cdot 2^{m-1}-1,
\end{eqnarray}
then $\rho\left(C_s((q_0)\right)=3$. 
\end{theorem}

\begin{proof} If $q_0=3$, then $m=1$, it is easy to check that Eq (\ref{x:bounds-odd}) holds for $s=3$. Since it was known that $\rho\left(C_s(3)\right)=3$ for any $s \ge 2$, Theorem \ref{x4:thm4-odd} holds for this case. From now on let us assume that $q_0 \ge 5$.

Denote by $N$ the number of $y \in \gf_{q}^*\setminus \square_{q_0}$ such that
\[y\left(y-\beta\right) \in \square_q, \quad \forall \beta \in \square_{q_0}.\]
We shall prove that $N>0$ if $s \ge 3$ is odd and $s$ satisfies Eq \eqref{x:bounds-odd}. This would imply $\rho\left(C_s(q_0)\right)=3$ thanks to Theorem \ref{r3s-odd}. 

Let 
\[\square_{q_0}=\left\{b_1,\cdots,b_m\right\}\]
and $\chi: \gf_{q}^* \to \left\{\pm 1\right\}$ be the quadratic character. We have 
\begin{eqnarray}
    N&=&\sum_{y \in \gf_{q}^* \setminus \square_{q_0}} \prod_{i=1}^m \frac{1}{2} \big\{1+\chi\left(y(y-b_i)\right)\big\} =\frac{1}{2^m} \left(A-B\right), \label{x:etN}
\end{eqnarray}
where
\begin{eqnarray}
    A&=&\sum_{y \in \gf_{q}^*} \prod_{i=1}^m \big\{1+\chi\left(y(y-b_i)\right)\big\}, \label{x:A} \\
    B&=&\sum_{y \in \square_{q_0}} \prod_{i=1}^m \big\{1+\chi\left(y(y-b_i)\right)\big\}. \label{x:B}\nonumber
\end{eqnarray}
We deal with $B$ first. Since $s$ is odd and $q=q_0^2$, $\square_{q_0} \subseteq \square_{q}$, so we have 
\begin{eqnarray*}
    B&=&\sum_{y \in \square_{q_0}} \prod_{i=1}^m \big\{1+\chi\left(y(y-b_i)\right)\big\}=\sum_{y \in \square_{q_0}} \prod_{\beta \in \square_{q_0}} \big\{1+\chi\left(y-\beta\right)\big\}.
\end{eqnarray*}
We claim that for any fixed $y \in \square_{q_0}$, there is $\beta \in \square_{q_0}$ such that $\chi(y-\beta)=-1$. The proof is simple: first, since $s$ is odd, $\chi$ can be considered as the quadratic character defined on $\gf_{q_0}$. Now suppose the claim is not true, then $\chi(y-\beta)=1$ for any $\beta \in \square_{q_0} \setminus \{y\}$, so as $q_0 \ge 5$ we have
\begin{eqnarray} \label{x:qcharsum} \sum_{x \in \gf_{q_0}} \chi\left(y-x^2\right)=1+2\sum_{\beta \in \square_{q_0} \setminus \{y\}} \chi\left(y-\beta\right)=1+2(m-1)=q_0-2 \ge 3.\end{eqnarray}
On the other hand, by Lemma \ref{x:lem-qua} we have
\[\left|\sum_{x \in \gf_{q_0}} \chi\left(y-x^2\right)\right|=1,\]
which clearly contradicts \eqref{x:qcharsum}. So the claim is proved, and hence $B=0$. 

Now we deal with the quantity $A$. Expanding the right hand side of Eq \eqref{x:A} and noting that $\chi$ is multiplicative, we have 
\begin{align}
A&=\sum_{y \in \gf_{q}^*} \left\{1+\sum_{\emptyset\neq I\subseteq [m]\,} \chi\left(y^{\#I}\prod_{i\in I}(y-b_i)\right)\right\} =q-1+A_1, \label{x:etA}
\end{align}
where
\begin{align} \label{x:A1-odd}
    A_1&=\sum_{\emptyset\neq I\subseteq [m]\,} \sum_{y \in \gf_{q}^*} \chi(y)^{\#I} \chi\left(\prod_{i\in I}(y-b_i)\right). 
\end{align}
Denote
\begin{align*}
    T&=\left|\sum_{y \in \gf_{q}^*} \chi(y)^{\#I} \chi\left(\prod_{i\in I}(y-b_i)\right)\right|.
\end{align*}
If $\#I$ is odd, then by Lemma \ref{weil} 
\begin{eqnarray} \label{x:T1-odd} T=\left|\sum_{y \in \gf_{q}} \chi\left(y\prod_{i\in I}(y-b_i)\right)\right| \le 1+\left(\#I+1-2\right) \sqrt{q}=1+\left(\#I-1\right) \sqrt{q}.\end{eqnarray}
If $\#I$ is even, then by Lemma \ref{weil} again 
\begin{eqnarray} \label{x:T2-odd} T&=&\left|\sum_{y \in \gf_{q}^*} \chi\left(\prod_{i\in I}(y-b_i)\right)\right| 
\le 1+ \left|\sum_{y \in \gf_{q}} \chi\left(\prod_{i\in I}(y-b_i)\right)\right| \nonumber \\
&\le &1+1+\left(\#I-2\right) \sqrt{q}=2+\left(\#I-2\right) \sqrt{q}.\end{eqnarray}
Plugging the estimates \eqref{x:T1-odd} and \eqref{x:T2-odd} into the right side of \eqref{x:A1-odd}, we obtain
\[|A_1| \leq \sum_{\substack{\emptyset \ne I \subseteq [m]\\
\#I \mbox{\small \, odd}}} \big\{1+\left(\#I-1\right)\sqrt{q}\big\}+\sum_{\substack{\emptyset \ne I \subseteq [m]\\
\#I \mbox{\small \, even}}} \big\{2+\left(\#I-2\right) \sqrt{q}\big\}.\]
From the elementary identities
\begin{align*}
    \sum_{\substack{\emptyset \ne I \subseteq [m]\\
\#I \mbox{\small \, odd}}}1&=\sum_{i \mbox{\small \, odd}} \binom{m}{i}=2^{m-1},\\
\sum_{\substack{\emptyset \ne I \subseteq [m]\\
\#I \mbox{\small \, even}}}1&=\sum_{0<i \mbox{\small \, even}} \binom{m}{i}=2^{m-1}-1,\\
\sum_{\emptyset \ne I \subseteq [m]}\#I&=\sum_{i=1}^m \binom{m}{i}i=m \cdot 2^{m-1},\end{align*}
we can obtain easily that 
\[|A_1| \le \left(3 \cdot 2^{m-1}-2 \right)+\sqrt{q}\big((m-3) \cdot 2^{m-1}+2 \big).\]
Finally, noting \eqref{x:etN}, \eqref{x:etA} and $B=0$,  we obtain 
\begin{eqnarray*}
2^m \cdot N=q-1+A_1-B \ge q-1-\left\{\left(3 \cdot 2^{m-1}-2 \right)+\sqrt{q}\big((m-3) \cdot 2^{m-1}+2 \big)\right\}.
\end{eqnarray*}
Now it is clear that if $s$ satisfies \eqref{x:bounds-odd}, then $2^m \cdot N>0$ and hence $N>0$ as desired. This completes the proof of Theorem \ref{x4:thm4-odd}. 
\end{proof}

It is interesting to compare Theorem \ref{x4:thm4-odd} with \cite[Theorem IV.4]{SHO23} and \cite[Theorem 4.8]{SLHO25} which state that if $q_0 \ge 5$ and $s \ge 3$ are both odd and satisfy 
	\begin{eqnarray} \label{x:shiine}
		q_0^s-q_0^{s/2}\left((m-2)2^{m}+2\right)>2^{m}-1,     
	\end{eqnarray}
then $\rho\left(C_s((q_0)\right)=3$. While the bound on $s$ from \eqref{x:bounds-odd} looks stronger than that from \eqref{x:shiine}, they are essentially of the same quality, as seen from the following lemma.

\begin{lemma} \label{x4:lem-odd} For any odd prime power $q_0$, denote by $s_*(q_0) \ge 3$ and $s'_{*}(q_0) \ge 3$ the least odd number satisfying \eqref{x:bounds-odd} and \eqref{x:shiine} respectively. If $q_0 \ge 13$, then \begin{itemize}
	\item[i).] $s_*(q_0)$ is an odd number in the range 
	\begin{eqnarray} \label{x:ss} 
		\frac{q_0 \log 2-5}{\log q_0}+2< s_*(q_0) < \frac{q_0\log 2-5 \log 2}{\log q_0}+4.
	\end{eqnarray}
	\item[ii).] $s'_*(q_0)-s_*(q_0) \in \{0,2\}$. 
	\end{itemize}
	\end{lemma}
\begin{proof}
i). For a fixed $q_0 \ge 13$, denote by $s_1,s_2$ the unique positive real numbers such that 
\begin{eqnarray} \label{x:s1s2-1} q_0^{s_1/2}=(m-3)2^{m-1}+2, \quad q_0^{s_2/2}=(m-3)2^{m-1}+3. \end{eqnarray}
We claim that
\begin{eqnarray} \label{x:s12-0} s_1<s_*(q_0)<s_2+2.\end{eqnarray}
The lower bound of $s_*(q_0)$ in \eqref{x:s12-0} is obvious. As for the upper bound, let $s':=\lceil s_2 \rceil \ge s_2$. By \eqref{x:s1s2-1} and nothing that $m \ge 6$ we have 
\begin{align*}
	q_0^{s'}-q_0^{s'/2}\left((m-3)2^{m-1}+2\right)&\ge q_0^{s_2}-q_0^{s_2/2}\left((m-3)2^{m-1}+2\right)\\
	&=q_0^{s_2/2}=(m-3)2^{m-1}+2>3 \cdot 2^{m-1}-1.\end{align*}
Since $s_*(q_0)$ is the least odd number satisfying Eq \eqref{x:bounds-odd}, we must have $s_*(q_0) \le \lceil s_2 \rceil+1 <s_2+2$, which is the desired upper bound. So the claim \eqref{x:s12-0} is proved. 

Now we estimate $s_1$ and $s_2$. From \eqref{x:s1s2-1} it is easy to check that 
\begin{eqnarray} \label{x:s1s2-2} 
	s_i&=&2+\frac{q_0\log 2-5\log 2}{\log q_0}+\frac{2 \log \left(1-\frac{7}{q_0}+\frac{\delta_i}{q_0\sqrt{2}^{q_0-5}}\right)}{\log q_0}, \quad \delta_1=2,\, \,\delta_2=3.
	\end{eqnarray}
It is easy to check, with the help of {\sc Magma} or {\sc Mathematica} that whenever $q_0 \ge 13$ the following hold:
\begin{eqnarray*}
	s_1-\frac{q_0 \log 2-5}{\log q_0}-2&=&\frac{1}{\log q_0}\left(5-5 \log 2+2 \log \left(1-\frac{7}{q_0}+\frac{2}{q_0\sqrt{2}^{q_0-5}}\right)\right)>0,\\
	s_2-\frac{q_0 \log 2-5\log 2}{\log q_0}-2&=&\frac{2}{\log q_0} \log \left(1-\frac{7}{q_0}+\frac{2}{q_0\sqrt{2}^{q_0-5}}\right)<0,
\end{eqnarray*}
that is, 
\begin{eqnarray*}
	\frac{q_0 \log 2-5}{\log q_0}+2<s_1<s_2<\frac{q_0 \log 2-5\log 2}{\log q_0}+2, \quad \forall q_0 \ge 13. 
\end{eqnarray*}
Plugging this into Eq \eqref{x:s12-0}, we obtain \eqref{x:ss} as desired. 

ii). Denote by $s_1',s_2'$ the unique positive real numbers such that 
\begin{eqnarray*} \label{x:s1s2-s1} q_0^{s_1/2}=(m-2)2^{m}+2, \quad q_0^{s_2/2}=(m-2)2^{m}+3. \end{eqnarray*}
By using the same argument as in i), we have 
\begin{eqnarray*} \label{x:s12-1} s_1'<s'_*(q_0)<s_2'+2, \end{eqnarray*}
and $s_1',s_2'$ can be computed as 
\begin{eqnarray} \label{x:s1s2-3} 
	s_i'&=&2+\frac{q_0\log 2-3\log 2}{\log q_0}+\frac{2 \log \left(1-\frac{5}{q_0}+\frac{\delta_i}{q_0\sqrt{2}^{q_0-3}}\right)}{\log q_0}, \quad \delta_1=2,\, \,\delta_2=3.
\end{eqnarray}
Thus we have
\[0 \le s'_*(q_0)-s_*(q_0)<s_2'+2-s_1. \]
Using the expressions of $s_1$ in \eqref{x:s1s2-2} and of $s_2'$ in \eqref{x:s1s2-3}, it is easy to check by {\sc Magma} that when $q_0 \ge 13$, 
\[s_2'-s_1=\frac{2}{\log q_0} \left(\log 2+\log \frac{1-\frac{5}{q_0}+\frac{3}{q_0\sqrt{2}^{q_0-3}}}{1-\frac{7}{q_0}+\frac{2}{q_0\sqrt{2}^{q_0-5}}}\right)<1, \]
so $0 \le s'_*(q_0)-s_*(q_0)<3$. Since both $s'_*(q_0)$ and $s_*(q_0)$ are odd, we have  $s'_*(q_0)-s_*(q_0) \in \{0,2\}$ as desired. 
\end{proof}

\begin{remark}
	In Table \ref{x:t2-odd} we list the values $s_*(q_0)$ and $s'_*(q_0)$ for all odd prime powers $5 \le q_0 \le 59$ which we computed by {\sc Magma}. The data in Table \ref{x:t2-odd} confirms Lemma \ref{x4:lem-odd}. 
	
	\begin{table}[h]
		\label{x:t2-odd}
		\caption{Comparison of $s_*(q_0)$ and $s'_*(q_0)$ which were derived from \eqref{x:bounds-odd} and \eqref{x:shiine} respectively. }
		\centering
		\begin{tabular}{|c|c|c|c|c|c|c|c|c|c|c|c|c|c|c|c|c|c|c|c|}
			\hline
			$q_0$ & 5 & 7 & 9&11&13&17&19&23&25&27&29&31&37&41&43&47&49&53&59\\ 
			\hline
			$s_*(q_0) $&3&3&3&3&5&5&5&7&7&7&7&9&9&9&9&11&11&11&13\\
			\hline
			$s'_*(q_0)$ &3&3&5&5&5&7&7&7&7&7&9&9&9&11&11&11&11&11&13\\
			\hline
		\end{tabular}
	\end{table}	
\end{remark}	
	
	We can summarize Corollaries \ref{x4:cor2-odd}, \ref{x4:cor3-odd}, Theorem \ref{x4:thm4-odd} and  Lemma \ref{x4:lem-odd} in the following proposition. For the sake of completeness, we also include the trivial case that $s=1$. For $q_0$ odd, let us define $s^*(q_0) \ge 3$ to be the largest odd number $s$ satisfying \eqref{x4:sodd}. 
	
	\begin{proposition} \label{x:prop-odd}
	Let $q_0$ be odd. The covering radius of $C_s(q_0)$ is given by 
	\begin{itemize}
		\item[i).] If $s=1$, then $ \rho\left(C_s(q_0)\right)=2$;
		\item[ii).] If $s\ge 2$ is even, then $\rho\left(C_s(q_0)\right)=3$;
		\item[iii).] If $s \ge 3$ is odd, and $s \le s^*(q_0)$, then $\rho\left(C_s(q_0)\right)=2$;
		\item[iv).] If $s \ge 3$ is odd, and  $s \ge s_*(q_0)$, then  $\rho\left(C_s(q_0)\right)=3$;
		\item[v).] If $s \ge 3$ is odd, and $s'|s$ with $\rho\left(C_{s'}(q_0)\right)=3$, then $\rho\left(C_s(q_0)\right)=3$.
	\end{itemize}
\end{proposition}

In Table \ref{x:t2-odd2} below, we list the values $s^*(q_0)$ and $s_*(q_0)$ for odd prime powers $q_0$ such that $3 \le q_0 \le 59$. We also list the missing odd numbers $s$ in the interval $(s^*(q_0),s_*(q_0))$ that $\rho\left(C_s(q_0)\right)$ is not yet determined.

\begin{table}[h]
	\label{x:t2-odd2}
	\caption{For each $q_0$, $s^*(q_0) \ge 3$ is the largest odd number satisfying \eqref{x4:sodd}, $s_*(q_0) \ge 3$ is the least odd number satisfying \eqref{x:bounds-odd}. }
	\centering
	\begin{tabular}{|c|c|c|c|}
		\hline
		$q_0$ & $s^*$ & $s_*$ & the set of missing odd numbers in $(s^*,s_*)$ \\ 
		\hline
		3 &  & 3& $\emptyset$\\
		\hline
		5 &  & 3& $\emptyset$\\
		\hline
		7 &  & 3& $\emptyset$\\
		\hline
		9 &  &3&  $\emptyset$\\
		\hline
		11 &   &3& $\emptyset$\\
		\hline
		13 &&5&  \{3\}\\
		\hline
		17 & & 5&  \{3\}\\
		\hline
		19 &  3 &5& $\emptyset$\\
		\hline
		23 &  3 &7& \{5\}\\
		\hline
		25 &  3 &7& \{5\}\\
		\hline
		27 &  3 &7& \{5\}\\
		\hline
		29 &  3 &7& \{5\}\\
		\hline
		31 &  3 &9& \{5,7\}\\
		\hline
		37 &  3 &9& \{5,7\}\\
		\hline
		41 &  3 &9& \{5,7\}\\
		\hline
		43 &   3&9& \{5,7\}\\
		\hline
		47 &   3&11& \{5,7,9\}\\
		\hline
		49 &   3&11& \{5,7,9\}\\
		\hline
		53 &   3&11& \{5,7,9\}\\
		\hline
		59 &  3 &13& \{5,7,9,11\}\\
		\hline
	\end{tabular}
\end{table}

The computations performed using {\sc Magma} agree perfectly with Proposition \ref{x:prop-odd}. For the missing odd $s$ in $(s^*(q_0),s_*(q_0))$, for example, we have $\rho(C_3(13))=3,\rho(C_3(17))=2, \rho(C_5(23))=3, \rho(C_5(25))=3,\rho(C_5(27))=3,\rho(C_5(29))=3$.

Since $\rho\left(C_s(q_0)\right)=\rho\left(\widetilde{C_s(q_0)}\right)$ when $q_0$ is odd, the following proposition holds.

	\begin{proposition} \label{x:prop-odd2}
	Let $q_0$ be odd. The covering radius of $\widetilde{C_s(q_0)}$ is given by 
	\begin{itemize}
		\item[i).] If $s=1$, then $ \rho\left(\widetilde{C_s(q_0)}\right)=2$;
		\item[ii).] If $s\ge 2$ is even, then $\rho\left(\widetilde{C_s(q_0)}\right)=3$;
		\item[iii).] If $s \ge 3$ is odd, and $s \le s^*(q_0)$, then $\rho\left(\widetilde{C_s(q_0)}\right)=2$;
		\item[iv).] If $s \ge 3$ is odd, and  $s \ge s_*(q_0)$, then  $\rho\left(\widetilde{C_s(q_0)}\right)=3$;
		\item[v).] If $s \ge 3$ is odd, and $s'|s$ with $\rho\left(\widetilde{C_{s'}(q_0)}\right)=3$, then $\rho\left(\widetilde{C_s(q_0)}\right)=3$.
	\end{itemize}
\end{proposition}

\section{Determining covering radius of $C_{s}(q_0)$: $q_0$ even}\label{sec6}

Let $ q_0 \geq 2 $ be even and $ s \geq 2 $. It is known by Corollary \ref{x3:cor} that $\rho(C_s(q_0)) \in \{2,3\}$. In this section we study explicit conditions to distinguish these two cases. 

We know by the proof of Theorem \ref{x2:thm} that $ \rho(C_s(q_0)) = 2 $ if and only if for any $ \alpha \in \mathbb{F}_q \setminus \mathbb{F}_{q_0} $ (note that the case $ \alpha \in \mathbb{F}_{q_0} $ is trivial), there exist $ (c_1, c_2) \in \mathbb{F}_{q_0}^{*2}$ and $ (x_1, x_2) \in H^2 $ such that 
\begin{eqnarray} \label{x5:r2}
\alpha = c_1 x_1 + c_2 x_2.
\end{eqnarray}
Taking the $ q $-th power on both sides of the above equation, we obtain
$$
\alpha = c_1 \overline{x}_1 + c_2 \overline{x}_2.
$$
Hence, $ x_1 $ satisfies the equation
$$
\left(c_1 x_1 + \alpha\right)(c_1 \overline{x}_1 + \alpha) = c_2^2,
$$
which simplifies to
\begin{eqnarray} \label{x4:eq1}
x_1^2 + \frac{\alpha^2 + c_1^2 + c_2^2}{\alpha c_1} \, \, x_1 + 1 = 0.
\end{eqnarray}
For the given $\alpha \in \gf_q \setminus \gf_{q_0}$ and $c_1,c_2 \in \gf_{q_0}^*$, if the above quadratic equation has a solution $x_1 \in H$, then by taking $x_2=\frac{\alpha+c_1x_1}{c_2}$, it is easy to see that $x_2 \in H$ and this pair $(x_1,x_2)\in H^2$ satisfies Eq (\ref{x5:r2}). So to study Eq (\ref{x5:r2}), we just need to study Eq (\ref{x4:eq1}). 

Since $ \frac{\alpha^2 + c_1^2 + c_2^2}{\alpha c_1} \neq 0 $ as $\alpha \in \gf_q \setminus \gf_{q_0}$ and $c_1,c_2 \in \gf_{q_0}^*$, by Lemma~\ref{rootinh2}, Eq (\ref{x4:eq1}) has a solution $ x_1 \in H $ if and only if 
$$
\mathrm{Tr}_{\gf_q/\gf_2}\left( \left(\frac{\alpha c_1}{(\alpha + c_1 + c_2)^2} \right)^2\right)=\mathrm{Tr}_{\gf_q/\gf_2}\left( \frac{\alpha c_1}{(\alpha + c_1 + c_2)^2} \right) = 1.
$$
Note that
\begin{align*}
	&\mathrm{Tr}_{\gf_q/\gf_2}\left(\frac{\alpha c_1}{(\alpha+c_1+c_2)^2}\right)\\
	=&\mathrm{Tr}_{\gf_q/\gf_2}\left(\frac{c_1(\alpha+c_1+c_2)}{(\alpha+c_1+c_2)^2}\right)+\mathrm{Tr}_{\gf_q/\gf_2}\left(\frac{ c_1(c_1+c_2)}{(\alpha+c_1+c_2)^2}\right)\\
	=&\mathrm{Tr}_{\gf_q/\gf_2}\left(\frac{c_1}{\alpha+c_1+c_2}\right)+\mathrm{Tr}_{\gf_q/\gf_2}\left(\frac{ c_1^2}{(\alpha+c_1+c_2)^2}\right)+\mathrm{Tr}_{\gf_q/\gf_2}\left(\frac{ c_1c_2}{(\alpha+c_1+c_2)^2}\right)\\
	=&\mathrm{Tr}_{\gf_q/\gf_2}\left(\frac{ c_1c_2}{(\alpha+c_1+c_2)^2}\right).
\end{align*}
The above derivation shows that the existence of a solution $ x_1 \in H $ to Eq (\ref{x4:eq1}) is equivalent to the condition:
$$
\mathrm{Tr}_{\gf_q/\gf_2}\left( \frac{c_1 c_2}{(\alpha + c_1 + c_2)^2} \right) = 1.
$$
This directly leads to the characterization of conditions under which $ \rho(C_s(q_0)) = 2 $. Specifically, for $ \rho(C_s(q_0)) = 2 $, this condition must hold for all $ \alpha \in \mathbb{F}_q \setminus \mathbb{F}_{q_0} $, while for $ \rho(C_s(q_0)) = 3 $, it must fail for at least one such $ \alpha $. We summarize these results in the following lemma.

\begin{lemma}\label{rho=3}
	Let $ q_0 $ be even and $ s \geq 2 $. Then $ \rho(C_s(q_0)) = 2 $ if and only if for every $ \alpha \in \mathbb{F}_q \setminus \mathbb{F}_{q_0} $, there exist $ c_1, c_2 \in \mathbb{F}_{q_0}^* $ such that 
	$$
	\mathrm{Tr}_{\gf_q/\gf_2}\left( \frac{c_1 c_2}{(\alpha + c_1 + c_2)^2} \right) = 1.
	$$
Otherwise $ \rho(C_s(q_0)) = 3 $ if and only if there exists $ \alpha \in \mathbb{F}_q \setminus \mathbb{F}_{q_0} $ such that \begin{equation}\label{thm2-6}
	\mathrm{Tr}_{\gf_q/\gf_2}\left( \frac{c_1 c_2}{(\alpha + c_1 + c_2)^2} \right) = 0, \quad \forall c_1,c_2 \in \gf_{q_0}^*.
	\end{equation}
\end{lemma}

Let us simplify Eq (\ref{thm2-6}) further in Lemma \ref{rho=3}. Suppose Eq (\ref{thm2-6}) holds for any $c_1,c_2 \in \gf_{q_0}^*$. Obvious it also holds when $c_1=0$ or $c_2=0$. Next for $c_1,c_2\in\gf_{q_0}$, define 
\begin{eqnarray} \label{x4:eq2} c_1c_2=a, \quad c_1+c_2=b,\end{eqnarray}
where $a,b \in \gf_{q_0}$. If $b=0$, then $c_1=c_2$, so in this case Eq (\ref{thm2-6}) reduces to $\mathrm{Tr}_{\gf_q/\gf_2}\left({c_1}/{\alpha}\right)=0$ for any $c_1 \in \gf_{q_0}$. Rewriting this as 
\begin{equation*}\label{thm2-7}
\mathrm{Tr}_{\gf_q/\gf_2}\left({c_1}/{\alpha}\right)=\mathrm{Tr}_{\gf_{q_0}/\gf_2}\left(c_1\mathrm{Tr}_{\gf_q/\gf_{q_0}}\left({1}/{\alpha}\right)\right)=0, \quad \forall c_1 \in \gf_{q_0},
\end{equation*}
we easily derive 
\begin{eqnarray} \label{x5:r31} \mathrm{Tr}_{\gf_q/\gf_{q_0}}\left({1}/{\alpha}\right)=0.\end{eqnarray}

If $b\neq 0$, then by Lemma \ref{rootinfq}, the quadratic equation $x^2+bx+a=0$ has two distinct roots $x=c_1,c_2$ in $\gf_{q_0}$ satisfying Eq (\ref{x4:eq2}) if and only if 
\begin{equation}\label{thm2-8}
 \mathrm{Tr}_{\gf_{q_0}/\gf_2}\left({a}/{b^2}\right)=0.
\end{equation}
Let $a=b^2t$. Then Eq (\ref{thm2-8}) implies that $\mathrm{Tr}_{\gf_{q_0}/\gf_2}(t)=0$, which means $t=\tau^2+\tau$ for some $\tau\in\gf_{q_0}$. Substituting $a=b^2t=b^2(\tau^2+\tau)$ into \eqref{thm2-6}, noting $c_1c_2=a,c_1+c_2=b$, we obtain
\[\mathrm{Tr}_{\gf_q/\gf_2}\left(\frac{b^2(\tau^2+\tau)}{(\alpha+b)^2}\right)=0, \quad \forall b \in \gf_{q_0}^*,\tau \in \gf_{q_0}.\]
This can be written as 
\[\mathrm{Tr}_{\gf_{q_0}/\gf_2}\left(\tau\,\,\mathrm{Tr}_{\gf_q/\gf_{q_0}}\left(\frac{b}{\alpha+b}+\frac{b^2}{(\alpha+b)^2}\right)\right)=0, \quad \forall b \in \gf_{q_0}^*,\tau \in \gf_{q_0}. \]
So we have 
\[\mathrm{Tr}_{\gf_q/\gf_{q_0}}\left(\frac{b}{\alpha+b}+\frac{b^2}{(\alpha+b)^2}\right)=0, \quad \forall b \in \gf_{q_0}^*,\]
from which we conclude easily that 
\[\mathrm{Tr}_{\gf_q/\gf_{q_0}}\left(\frac{b}{\alpha+b}\right)\in\{0,1\}, \quad \forall b \in \gf_{q_0}^*.\]

Combining this and Eq \eqref{x5:r31}, we find that Eq (\ref{thm2-6}) in Lemma~\ref{rho=3} is equivalent to 
\begin{enumerate}
	\item $\mathrm{Tr}_{\gf_q/\gf_{q_0}}\left( \frac{1}{\alpha} \right) = 0$, and 
	\item $\mathrm{Tr}_{\gf_q/\gf_{q_0}}\left( \frac{b}{\alpha + b} \right) \in \{0, 1\}$ for all $b \in \mathbb{F}_{q_0}^*$.
\end{enumerate}
By substitution $\alpha\to\frac{1}{\alpha}$, we obtain the following simplified characterization.

\begin{theorem}
	\label{rho3simplified}
	Let $q_0 \ge 2$ be even and $s \geq 2$. Then $\rho(C_s(q_0)) = 3$ if and only if there exists $\alpha \in \mathbb{F}_q \setminus \mathbb{F}_{q_0}$ such that
	\begin{enumerate}
		\item $\mathrm{Tr}_{\gf_q/\gf_{q_0}}(\alpha) = 0$, and 
		\item $\mathrm{Tr}_{\gf_q/\gf_{q_0}}\left( \frac{1}{1 + b\alpha} \right) \in \{0, 1\}$ for all $b \in \mathbb{F}_{q_0}^*$.
	\end{enumerate}
Otherwise $\rho(C_s(q_0))=2$.     
\end{theorem}

As simple consequences of Theorem~\ref{rho3simplified}, we collect below several results about $\rho(C_s(q_0))$ for some special $q_0$ and $s$.

\begin{corollary} \label{x4:cor2} Let $q_0$ be even and $s \ge 2$ be a positive integer. If $s=2$ or if $s \ge 3$ is odd and $s \le \frac{q_0}{2}$, then $\rho(C_s(q_0))=2$. 
\end{corollary}

\begin{proof}
If $s=2$ then $q=q_0^2$, for any $\alpha \in \gf_q$, 
\[\mathrm{Tr}_{\gf_q/\gf_{q_0}}(\alpha)=\alpha+\alpha^{q_0}=0 \quad \Longrightarrow \quad \alpha \in \gf_{q_0}.\]
So there is no $\alpha \in \gf_{q}\setminus \gf_{q_0}$ such that $\mathrm{Tr}_{q_0}^q(\alpha)=0$, by Theorem \ref{rho3simplified}, $\rho(C_s(q_0))=2$. 

Now suppose $s$ is odd and $s \le q_0/2$. For any $\alpha \in \gf_{q} \setminus \gf_{q_0}$, let $t$ be the least positive integer such that $\alpha^{q_0^t}=\alpha$. It is easy to see that $t|s$ and  $\alpha,\ldots, \alpha^{q_0^{t-1}}$ are distinct elements in $\gf_{q_0^t} \subseteq \gf_{q}$. So the rational function $f(x) \in \gf_q(x)$ defined by 
\[f(x):=\sum_{i=0}^{s-1}\frac{1}{1+x\alpha^{q_0^i}}\]
satisfies 
\[f(x)=\frac{s}{t} \left(\sum_{i=0}^{t-1}\frac{1}{1+x\alpha^{q_0^i}}\right)=\sum_{i=0}^{t-1}\frac{1}{1+x\alpha^{q_0^i}}.\]
Let us define 
\[g(x)=\sum_{i=0}^{t-1}\frac{1}{1+x\alpha^{q_0^i}}.\]
Then $g(x)$ is a non-constant rational function defined over $\gf_{q}$. We can write 
\[g(x)=\frac{N(x)}{D(x)},\]
where $D(x),N(x) \in \gf_{q}[x]$ are given by 
\[D(x)=\prod_{i=0}^{t-1} \left(1+x \alpha^{q_0^i}\right), \quad N(x)=D(x) \left(\sum_{i=0}^{t-1} \frac{1}{1+x\alpha^{q_0^i}}\right).\]
It is easy to see that $\deg D=t$ and $\deg N \le t-1$, and $N\ne 0$. 

For any $\epsilon$, denote by $N_{\epsilon}$ the number of solutions $x \in \gf_{q_0}$ such that $f(x)=g(x)=\epsilon$, which is equivalent to the equation
\[N(x)+\epsilon D(x)=0.\]
Since $N(x),D(x) \in \gf_q[x]$ are not zero polynomials with $\deg D=t, \, \deg N \le t-1$, we have 
\[N_0+N_1 \le t-1+t=2t-1 \le 2s-1 \le q_0-1,\]
which implies that there is always $b \in \gf_{q_0}$ such that 
\[\mathrm{Tr}_{\gf_q/\gf_{q_0}}\left(\frac{1}{1+b\alpha }\right) =\frac{1}{1+b\alpha }+\frac{1}{1+b\alpha^{q_0} }+\cdots+\frac{1}{1+b\alpha^{q_0^{s-1}} } =f(b) \not \in \{0,1\}.\]
Obvioulsy $b \ne 0$. Thus by Theorem \ref{rho3simplified}, $\rho(C_s(q_0))=2$. 
\end{proof}

\begin{corollary} \label{x4:cor3}
	Let $q_0$ be even. 
    \begin{itemize}
        \item[1).] If $s\geq 4$ is even, then $\rho(C_s(q_0))=3$.
        \item[2).] If $s$ is odd, and there is $s'|s$ such that $\rho(C_{s'}(q_0))=3$, then $\rho(C_s(q_0))=3$. 
    \end{itemize}
    \end{corollary}
\begin{proof}
1). Since $s \ge 4$ is even, we have $q=q_0^s=q_1^2$, where $q_1=q_0^{\frac{s}{2}}>q_0$. Take any $\alpha\in \gf_{q_1}\setminus\gf_{q_0}$. Then
\[\mathrm{Tr}_{\gf_q/\gf_{q_0}}\left(\alpha\right)=\mathrm{Tr}_{\gf_{q_1}/\gf_{q_0}}\left(\alpha\mathrm{Tr}_{\gf_q/\gf_{q_1}}(1)\right)=0,\]
and
\[\mathrm{Tr}_{\gf_q/\gf_{q_0}}\left(\frac{1}{1+b\alpha}\right)=\mathrm{Tr}_{\gf_{q_1}/\gf_{q_0}}\left(\frac{1}{1+b\alpha}\mathrm{Tr}_{\gf_q/\gf_{q_1}}(1)\right)=0, \quad \forall b \in \gf_{q_0}^*.\]
By Theorem \ref{rho3simplified}, $\rho(C_s(q_0))=3$.

2). By Theorem \ref{rho3simplified}, since $\rho(C_{s'}(q_0))=3$, then there exists $\alpha\in\gf_{q_0^{s'}}\setminus\gf_{q_0}$, such that $\mathrm{Tr}_{\gf_{q_0^{s'}}/\gf_{q_0}}\left(\alpha\right)=0$, and $\mathrm{Tr}_{\gf_{q_0^{s'}}/\gf_{q_0}}\left(\frac{b}{\alpha+b}\right)\in\{0,1\}$ for all $b\in\gf_{q_0}^*$. Since $s'|s$ and $s$ is odd, we have $\gf_{q_0^{s'}}\subseteq\gf_{q_0^s}$, and $\frac{s}{s'}$ is odd. So
\[\mathrm{Tr}_{\gf_{q_0^{s}}/\gf_{q_0}}\left(\alpha\right)=\mathrm{Tr}_{\gf_{q_0^{s'}}/\gf_{q_0}}\left(\alpha\mathrm{Tr}_{\gf_{q_0^{s}}/\gf_{q_0^{s'}}}(1)\right)=0\]
and
\[\mathrm{Tr}_{\gf_{q_0^{s}}/\gf_{q_0}}\left(\frac{1}{1+b\alpha}\right)=\mathrm{Tr}_{\gf_{q_0^{s'}}/\gf_{q_0}}\left(\frac{1}{1+b\alpha}\mathrm{Tr}_{\gf_{q_0^{s}}/\gf_{q_0^{s'}}}(1)\right)\in\{0,1\} \quad \forall b \in \gf_{q_0}^*.\]
Thus by Theorem \ref{rho3simplified}, we have $\rho(C_{s}(q_0))=3$. 
\end{proof}

Finally, for any $q_0 \ge 4$ even and $s$ odd, we prove that if $s$ is sufficiently large, then $\rho(C_s(q_0))=3$.



\begin{theorem} \label{x4:thm4} Let $q_0$ be even. If $s \ge 3$ is odd and satisfies
	\begin{eqnarray} \label{x:bounds-even}
		q_0^s-q_0^{s/2} \left(\frac{q_0}{2}\right)^{q_0-1} \Big(2q_0^2-7q_0+3\Big)>\left(\frac{q_0}{2}\right)^{q_0-1}\Big(4q_0-2\Big)-q_0^2,
	\end{eqnarray}
	then $\rho\left(C_s(q_0)\right)=3$. 
	\end{theorem}

\begin{proof} First, when $q_0=2$, it is well-known that $\rho(C_2(2))=2$ and $\rho(C_s(2))=3$ when $s \ge 3$ (see \cite{D85}) and it satisfies Theorem \ref{x4:thm4}, from now on let us assume that $q_0 \ge 4$ is even and $s \ge 3$ is odd. 

Since $q=q_0^s$, for simplicity define $\gf:=\gf_{q} \setminus \gf_{q_0}$ and denote by $N$ the number of $y \in \gf$ such that $\mathrm{Tr}_{\gf_{q}/\gf_{q_0}}(y) = 0$ and $\mathrm{Tr}_{\gf_{q}/\gf_{q_0}}\left( \frac{1}{1 + by} \right) \in \{0, 1\}$ for all $b \in \mathbb{F}_{q_0}^*$. We shall prove that $N>0$ if $s \ge 2q_0+1$, which implies that $\rho(C_s(q_0))=3$ thanks to Theorem \ref{rho3simplified}. 

Now we compute the quantity $N$. Set $\gf_{q_0}^*=\{b_1,b_2,\cdots,b_{r}\}$ where $r=q_0-1$. For any $\vec{e}=(e_1,\cdots,e_r) \in \{0,1\}^r$, denote by $N_{\vec{e}}$ the number of $y \in \gf$ such that 
\[\mathrm{Tr}_{\gf_{q}/\gf_{q_0}}(y) = 0 \mbox{ and } \mathrm{Tr}_{\gf_{q}/\gf_{q_0}}\left( \frac{1}{1 + b_iy} \right)=e_i, \quad \forall 1 \le i \le r.  \]
By the orthogonal relations of character sums over $\gf_q$, we have 
\[N_{\vec{e}}=\sum_{y\in\gf}\left(\frac{1}{q_0}\sum_{a\in\gf_{q_0}}(-1)^{\mathrm{Tr}_{\gf_{q_0}/\gf_{2}}\left(a\mathrm{Tr}_{\gf_{q}/\gf_{q_0}}(y)\right)}\right)\prod_{i=1}^{r}
	\frac{1}{q_0}\sum_{a_i\in\gf_{q_0}}(-1)^{\mathrm{Tr}_{\gf_{q_0}/\gf_{2}}\left(a_i\left(\mathrm{Tr}_{\gf_{q}/\gf_{q_0}}\left(\frac{1}{1+b_iy}\right)+e_i\right)\right)}.\]
Using the standard additive character $\psi$ on $\gf_{q}$ given by $\psi(x)=(-1)^{\mathrm{Tr}_{\gf_q/\gf_2}(x)}$, we can write  
\begin{align*}
	N_{\vec{e}}=&\frac{1}{q_0^{r+1}}\sum_{y\in\gf}\left(1+\sum_{a\in\gf_{q_0}^*}\psi(ay)\right)\prod_{i=1}^r\left\{1+\sum_{a_i\in\gf_{q_0}^*}\psi\left(a_i\left(\frac{1}{1+b_iy}+e_i\right)\right)\right\}.
\end{align*}
The product on the right can be expanded as 
\begin{align*}
1+\sum_{\emptyset\neq I\subseteq [r]\,} \sum_{\substack{a_i\in\gf_{q_0}^*\\
\forall i \in I}}\psi\left(\sum_{i\in I}a_i\left(\frac{1}{1+b_iy}+e_i\right)\right)=1+\sum_{\emptyset\neq I\subseteq [r]}\sum_{\substack{a_i\in\gf_{q_0}^*\\
\forall i \in I}} \psi\left(\sum_{i\in I}a_ie_i\right)\psi\left(\sum_{i\in I}\frac{a_i}{1+b_iy}\right).
\end{align*}
So we have 
\begin{eqnarray*}
	N&=&\sum_{\vec{e} \in \{0,1\}^r} N_{\vec{e}} \\
    &=&\frac{1}{q_0^{r+1}}\sum_{y\in\gf}\left(1+\sum_{a\in\gf_{q_0}^*}\psi(ay)\right) \sum_{\vec{e} \in\{0,1\}^r} \left(1+\sum_{\emptyset\neq I\subseteq [r]}
    \sum_{\substack{a_i\in\gf_{q_0}^*\\
\forall i \in I}} \psi\left(\sum_{i\in I}a_ie_i\right)\psi\left(\sum_{i\in I}\frac{a_i}{1+b_iy}\right)\right). 
\end{eqnarray*}
The sum over $\vec{e} \in\{0,1\}^r$ on the right is
\begin{eqnarray*}
2^r+\sum_{\emptyset\neq I\subseteq [r]}2^{r-\#I}
    \sum_{\substack{a_i\in\gf_{q_0}^*\\
\forall i \in I}} \sum_{\substack{e_i \in \{0,1\}\\
\forall i \in I}}\psi\left(\sum_{i\in I}a_ie_i\right)\psi\left(\sum_{i\in I}\frac{a_i}{1+b_iy}\right)\\
=2^r+\sum_{\emptyset\neq I\subseteq [r]}2^{r-\#I}
    \sum_{\substack{a_i\in\gf_{q_0}^*\\
\forall i \in I}} \left\{\prod_{i \in I}\left(1+\psi(a_i)\right)\right\} \psi\left(\sum_{i\in I}\frac{a_i}{1+b_iy}\right).
\end{eqnarray*}
Denote
\[\mathbb{H}_0:=\left\{a \in \gf_{q_0}^*: \mathrm{Tr}_{\gf_{q_0}/\gf_{2}}(a)=0\right\}.\]
It is easy to see that $\#\mathbb{H}_0=q_0/2-1$ and 
\[\prod_{i \in I}\left(1+\psi(a_i)\right)=\left\{\begin{array}{cl}
2^{\#I}     & \mbox{ if } a_i \in \mathbb{H}_0 \, \,\,\,\forall \, i, \\
0     & \mbox{ otherwise }. 
\end{array}\right.\]
So we have 
\begin{eqnarray}
N&=& \frac{2^r}{q_0^{r+1}}\sum_{y\in\gf}\left(1+\sum_{a\in\gf_{q_0}^*}\psi(ay)\right)  \left(1+\sum_{\emptyset\neq I\subseteq [r]}
    \sum_{\substack{a_i\in\mathbb{H}_0\\
\forall i \in I}} \psi\left(\sum_{i\in I}\frac{a_i}{1+b_iy}\right)\right)\nonumber\\
&=&\frac{2^r}{q_0^{r+1}} \Big(I_0+I_1+I_2+I_3\Big) \label{x4:vN}
\end{eqnarray}
where 
\begin{eqnarray}
I_0&=&\sum_{y \in \gf} 1=q-q_0 \label{x4:et0},\\
I_1&=&\sum_{a \in \gf_{q_0}^*} \sum_{y \in \gf} \psi(ay), \nonumber\\
I_2&=&\sum_{\emptyset\neq I\subseteq [r]}
    \sum_{\substack{a_i\in\mathbb{H}_0\\
\forall i \in I}} \sum_{y \in \gf} \psi\left(\sum_{i\in I}\frac{a_i}{1+b_iy}\right),\nonumber\\
I_3&=&\sum_{a\in\gf_{q_0}^*}\sum_{\emptyset\neq I\subseteq [r]}
    \sum_{\substack{a_i\in\mathbb{H}_0\\
\forall i \in I}} \sum_{y \in \gf} \psi\left(ay+\sum_{i\in I}\frac{a_i}{1+b_iy}\right). \nonumber
\end{eqnarray}
Since $s$ is odd, it is easy to see that 
\begin{eqnarray}I_1&=&\sum_{a \in \gf_{q_0}^*} \left(\sum_{y \in \gf_{q}} \psi(ay)-\sum_{y \in \gf_{q_0}} \psi(ay)\right)=0. \label{x4:et1} \end{eqnarray}
We can use Lemma \ref{rational} to give an estimate of $I_2$ and $I_3$ respectively. In particular, for $I_2$, for each $\emptyset \neq I \subseteq [r]$ and fixed $a_i \in \mathbb{H}_0$ where $i \in I$, for the rational function 
\[f(y)=\sum_{i \in I} \frac{a_i}{1+b_iy} \in \gf_q(y),\]
this corresponds to $M=0$, $L=2\#I$ in Lemma \ref{rational}, and consequently
\[\left|\sum_{y\in\gf_{q}\setminus\left\{b_i^{-1}: i \in I\right\}}
\psi\left(\sum_{i\in I}\frac{a_i}{1+b_iy}\right)\right|\leq 1+\left(2\#I-2\right)\sqrt{q}.\]
Hence we have
\begin{eqnarray*} \left|\sum_{y\in\gf}
\psi\left(\sum_{i\in I}\frac{a_i}{1+b_iy}\right)\right| &\le &q_0-\#I+ \left|\sum_{y\in\gf_{q}\setminus\left\{b_i^{-1}: i \in I\right\}}
\psi\left(\sum_{i\in I}\frac{a_i}{1+b_iy}\right)\right| \nonumber\\
&\le &q_0-\#I+1+2(\#I-1) \sqrt{q}. 
\end{eqnarray*}
Therefore
\begin{eqnarray}
    |I_2|&\le& \sum_{\emptyset\neq I\subseteq [r]}
    \sum_{\substack{a_i\in\mathbb{H}_0\\
\forall i \in I}} \left|\sum_{y \in \gf} \psi\left(\sum_{i\in I}\frac{a_i}{1+b_iy}\right)\right| \nonumber\\
&\le& \sum_{\emptyset\neq I\subseteq [r]}
    \sum_{\substack{a_i\in\mathbb{H}_0\\
\forall i \in I}} \Big(2(\#I-1) \sqrt{q}+q_0-\#I+1\Big)\nonumber\\
&=& \sum_{i=1}^r \binom{r}{i} \left(\frac{q_0}{2}-1\right)^i \Big(2(i-1) \sqrt{q}+q_0-i+1\Big). \label{x4:et2}
\end{eqnarray}
As for $I_3$, for each $\emptyset \neq I \subseteq [r]$, fixed $a \in \gf_{q_0}^*$ and fixed $a_i \in \mathbb{H}_0$ where $i \in I$, for the rational function 
\[g(y)=ay+\sum_{i \in I} \frac{a_i}{1+b_iy} \in \gf_q(y),\]
this corresponds to $M=1$, $L=2\#I$ in Lemma \ref{rational}, and consequently
\[\left|\sum_{y\in\gf_{q}\setminus\left\{b_i^{-1}: i \in I\right\}}
\psi\left(ay+\sum_{i\in I}\frac{a_i}{1+b_iy}\right)\right|\leq 1+\left(2\#I-1\right)\sqrt{q}.\]
Hence we have
\begin{eqnarray*} \left|\sum_{y\in\gf}
\psi\left(ay+\sum_{i\in I}\frac{a_i}{1+b_iy}\right)\right| &\le &q_0-\#I+ \left|\sum_{y\in\gf_{q}\setminus\left\{b_i^{-1}: i \in I\right\}}
\psi\left(ay+\sum_{i\in I}\frac{a_i}{1+b_iy}\right)\right| \nonumber\\
&\le &q_0-\#I+1+(2\#I-1) \sqrt{q}. 
\end{eqnarray*}
Therefore
\begin{eqnarray}
    |I_3|&\le& \sum_{a \in \gf_{q_0}^*}\sum_{\emptyset\neq I\subseteq [r]}
    \sum_{\substack{a_i\in\mathbb{H}_0\\
\forall i \in I}} \left|\sum_{y \in \gf} \psi\left(ay+\sum_{i\in I}\frac{a_i}{1+b_iy}\right)\right| \nonumber\\
&\le& (q_0-1)\sum_{\emptyset\neq I\subseteq [r]}
    \sum_{\substack{a_i\in\mathbb{H}_0\\
\forall i \in I}} \Big((2\#I-1) \sqrt{q}+q_0-\#I+1\Big)\nonumber\\
&=& (q_0-1)\sum_{i=1}^r \binom{r}{i} \left(\frac{q_0}{2}-1\right)^i \Big((2i-1) \sqrt{q}+q_0-i+1\Big). \label{x4:et3}
\end{eqnarray}

By \eqref{x4:et2} and \eqref{x4:et3} and a routine computation we have 
\begin{eqnarray}
    |I_2|+|I_3|&\le& \sum_{i=1}^r \binom{r}{i} \left(\frac{q_0}{2}-1\right)^i \Big(2(i-1) \sqrt{q}+q_0-i+1+(q_0-1)\left\{(2i-1) \sqrt{q}+q_0-i+1\right\}\Big) \nonumber\\
    &=&A_1\sqrt{q}+A_2, \label{x4:I23}
\end{eqnarray}
where
\begin{eqnarray*}
    A_1&=&\left(\frac{q_0}{2}\right)^{q_0-1} \Big(2q_0^2-7q_0+3\Big),\\
    A_2&=&\left(\frac{q_0}{2}\right)^{q_0-1}\Big(4q_0-2\Big)-q_0^2-q_0. 
\end{eqnarray*}

Combining \eqref{x4:vN}--\eqref{x4:et1} and \eqref{x4:I23}, we obtain
\begin{eqnarray*}
\frac{q_0^{r+1}}{2^r}N&=& I_0+I_1+I_2+I_3 \ge q-q_0-|I_2|-|I_3|\\
&\ge &q-q_0-A_1\sqrt{q}-A_2>0, 
\end{eqnarray*}
which is exactly Eq \eqref{x:bounds-even}. This completes the proof of Theorem \ref{x4:thm4}. 
\end{proof}

For $q_0$ even, denote by $s_*(q_0) \ge 3$ the least odd number satisfying \eqref{x:bounds-even}. We have the following result.  

\begin{lemma} \label{x4:lem-even} Let $q_0$ be even. Then 
	\begin{itemize}
		\item[i).] $s_*(2)=3$ and $s_*(4)= 7$; 
		\item[ii).] If $q_0 \ge 8$, then $s_*(q_0)$ is an odd number in the range 
		\begin{eqnarray*} \label{x:ss-even} 
			2q_0+2-\frac{2(q_0-1)\log 2}{\log q_0}< s_*(q_0) <2q_0+4-\frac{2(q_0-2)\log 2}{\log q_0}.
		\end{eqnarray*}
	\end{itemize}
\end{lemma}

\begin{proof}
i) can be verified by {\sc Magma} easily. The proof of ii) is elementary and very similar to that of \eqref{x:ss} in Lemma \ref{x4:lem-odd}. To be more precise, for fixed $q_0$, denote by $s_1,s_2$ the unique positive real numbers such that 
\begin{eqnarray*} \label{x:s1s2-even} q_0^{s_i/2}=\left(\frac{q_0}{2}\right)^{q_0-1}\left(2q_0^2-7q_0+3\right)+\delta_i, \quad \delta_1=0,\,\, \delta_1=1. \end{eqnarray*}
Very similar to the proof of \eqref{x:ss}, if $q_0 \ge 8$ we can obtain $s_1<s_*(q_0)<s_2+2$ and 
\[s_1>2q_0+2-\frac{2(q_0-1)\log 2}{\log q_0}, \quad s_2<2q_0+2-\frac{2(q_0-2)\log 2}{\log q_0}.\] 
This completes the proof of Lemma \ref{x4:lem-even}. 
\end{proof}

Let us also define $s^*(q_0) \ge 3$ to be the largest odd number satisfying $s^*(q_0) \le q_0/2$. Obviously $s^*(q_0)=q_0/2-1$ when $q_0 \ge 8$. We can summarize Corollaries \ref{x4:cor2}, \ref{x4:cor3}, and Theorem \ref{x4:thm4} in the following proposition. For the sake of completeness, we also include the trivial case that $s=1$. 

\begin{proposition} \label{x:prop-even}
	Let $q_0$ be even. The covering radius of $C_s(q_0)$ is given by 
	\begin{itemize}
		\item[i).] If $s=1$, then $ \rho\left(C_s(q_0)\right)=1$;
		\item[ii).] If $s=2$, then $\rho\left(C_s(q_0)\right)=2$;
		\item[iii).] If $s \ge 4$ is even, then $\rho\left(C_s(q_0)\right)=3$;
		\item[iv).] If $s \ge 3$ is odd, and  $s \le s^*(q_0)$, then  $\rho\left(C_s(q_0)\right)=2$;
		\item[v).] If $s \ge 3$ is odd, and  $s \ge s_*(q_0)$, then $\rho\left(C_s(q_0)\right)=3$;
		\item[vi).] If $s \ge 3$ is odd, and $s'|s$ with $\rho\left(C_{s'}(q_0)\right)=3$, then $\rho\left(C_s(q_0)\right)=3$.
		\end{itemize}
	\end{proposition}
	
	In Table \ref{x:t2-even} below, we list the values $s^*(q_0)$ and $s_*(q_0)$ for $q_0=2^t$ where $1 \le t \le 7$. We also list the missing odd numbers $s$ in the interval $(s^*(q_0),s_*(q_0))$ where the value $\rho\left(C_s(q_0)\right)$ is not yet determined. 
	
	\begin{table}[h]
		\label{x:t2-even}
		\caption{For each $q_0=2^t$, $s^*(q_0) \ge 3$ is the largest odd number $s$ such that $s \le q_0/2$, $s_*(q_0) \ge 3$ is the least odd number $s$ satisfying \eqref{x:bounds-even}. }
		\centering
		\begin{tabular}{|c|c|c|c|}
			\hline
			$q_0=2^t$ & $s^*$ & $s_*$ & the set of missing odd numbers in $(s^*,s_*)$ \\ 
			\hline
			2&  & 3&$\emptyset$ \\
			\hline
			4&  & 7& \{3,5\}\\
			\hline
			8 &3  &15& \{5,7,9,11,13\}\\
			\hline
			16 & 7 & 27& odd numbers from 9 to 25\\
			\hline
			32 & 15 &55& odd numbers from 17 to 53\\
			\hline
			64 & 31 &111& odd numbers from 33 to 109\\
			\hline
			128 &63&223&  odd numbers from 65 to 221\\
			\hline
		\end{tabular}
	\end{table}
	
The computations performed using {\sc Magma} agree perfectly with Proposition \ref{x:prop-even}. For the missing odd $s$ in $(s^*(q_0),s_*(q_0))$, for example, for $q_0=4$, we have $\rho(C_3(4))=2$ and $\rho(C_5(4))=3$; and for $q_0=8$, $\rho(C_5(8))=\rho(C_7(8))=2$.

\section{Discussion: quasi-perfect and maximal codes}\label{sec7}

Recall that quasi-perfect and maximal codes are  characterized by the minimum distance and the covering radius of the codes as were see in Eqs \eqref{x1:qp} and \eqref{x1:mc}. Given the detailed descriptions of the minimum distance and the covering radius in Theorems \ref{x3:thm1}-\ref{x3:thm2} and Propositions  \ref{x:prop-odd2}-\ref{x:prop-even}, we can identify quasi-perfect and maximal codes from $C_s(q_0)$ for even $q_0$ and from $\widetilde{C_s(q_0)}$ for odd $q_0$. 
	
For the sake of completeness, we also include the cases $s=1$ and $q_0=2,3$, which are either trivial or well-known. 

\subsection{$C_s(q_0)$ for even $q_0$}

{\bf Case $s=1$}: for $q_0$ even, $C_1(q_0)$ is a linear code with parameters $\left[q_0+1,q_0-1,3\right]_{q_0}$, $\rho(C_1(q_0))=1$, so $C_1(q_0)$ is both perfect and maximal. Actually $C_1(q_0)$ is equivalent to the Hamming code $\mathcal{H}_{q_0,2}$ (see \cite{HP03}) which is consistent with the classification of perfect codes. 

Now we consider the case that $s \ge 2$. 

{\bf Case $q_0=2$}: If $s=2$: $C_2(2)$ is a linear code with parameters $\left[5,1,5\right]_{q_0}$ with covering radius $\rho(C_1(q_0))=2$, so $C_2(2)$ is both perfect and maximal. The parameters of $C_2(2)$ are also consistent with the classification of perfect codes. 

If $s \ge 4$ even, $C_s(2)$ is a linear code with parameters $\left[2^s+1,2^s+1-2s,5\right]_{q_0}$ with covering radius $\rho(C_s(q_0))=3$, so $C_s(2)$ is both quasi-perfect and maximal. This result is well-known. 

{\bf Case $q_0 \ge 4$ is even}: $C_s(q_0)$ is a linear code with parameters $\left[2^s+1,2^s+1-2s\right]_{q_0}$. The minimum distance satisfies $d(C_s(q_0)) \in \{3,4\}$ and the covering radius satisfies $\rho(C_s(q_0))\in \{2,3\}$, so if $\rho(C_s(q_0))=2$, then $C_s(q_0)$ is quasi-perfect. In particular, 
\begin{itemize}
    \item[(1)] $s=2$: $d(C_s(q_0))=4, \rho(C_s(q_0))=2$, so $C_2(q_0)$ is both quasi-perfect and maximal; 
    \item[(2)] $s \ge 4$ even: $d(C_s(q_0))=4, \rho(C_s(q_0))=3$, so $C_2(q_0)$ is maximal;
    \item[(3)] $s \ge 3$ odd and $s \le \frac{q_0}{2}$: $d(C_s(q_0))=3, \rho(C_s(q_0))=2$, so $C_s(q_0)$ is both quasi-perfect and maximal. 
\end{itemize} 
We summarize the above result below (for simplicity, we focus on the case $q_0 \ge 4$ is even): 
\begin{theorem} \label{x5:cond}
    Let $q_0 \ge 4$ be even. 
    \begin{itemize}
    \item[(0)] If $s=1$, then $C_s(q_0)$ is both perfect and maximal; 
    \item[(1)] If $s=2$ or if $s$ is odd and $3 \le s \le q_0/2$, then $C_s(q_0)$ is both quasi-perfect and maximal; 
    \item[(2)] If $s \ge 4$ is even, then $C_s(q_0)$ is maximal. 
\end{itemize} 
\end{theorem}

The above results are summarized in the following table.

	\begin{table}[h]
	\label{t1}
        	\caption{ Parameters of $C_s(q_0)$ in even characteristic}
	\centering
	\begin{tabular}{|c|c|c|c|c|c|}
		\hline
		 $q_0$ even& $s$ & $d(C_s(q_0))$ & $\rho(C_s(q_0))$ & perfect/quasi-perfect & maximal \\ 
		\hline
		  any $q_0$ even& $1$ & $3$ & $1$ & perfect & $\checkmark$\\
		\hline
			$q_0=2$ & $2$ & $5$ & $2$ & perfect & $\checkmark$\\
		\hline
				$q_0=2$ & $s \ge 4$ and even & $5$ & $3$ & quasi-perfect & $\checkmark$\\
		\hline
					$q_0 \ge 4$ & $2$ & $4$ & $2$ & quasi-perfect & $\checkmark$\\
		\hline
				$q_0 \ge 4$ & $3\leq s \leq \frac{q_0}{2}$,~$s$~odd & $3$ & $2$ & quasi-perfect & $\checkmark$\\
		\hline
$q_0 \ge 4$ & $s\geq 4$ and~even & $4$ & $3$ & - & $\checkmark$\\
		\hline
	\end{tabular}
\end{table}

\subsection{$\widetilde{C_s(q_0)}$ for odd $q_0$}

{\bf Case $q_0=3$}: If $s=1$, $\widetilde{C_s(q_0)}$ is trivial with parameters $[2,0]$; if $s \ge 2$, $\widetilde{C_s(3)}$ is a linear code with parameters $\left[\frac{3^s+1}{2},\frac{3^s+1}{2}-2s,5\right]_3$, with covering radius $\rho\left(\widetilde{C_s(3)}\right)=3$, so $\widetilde{C_s(3)}$ is both quasi-perfec and maximal for any $s \ge 2$. This result is well-known. 

Now we assume that $q_0 \ge 5$. 

{\bf Case $s=1$}: $\widetilde{C_1(q_0)}$ is a linear code with parameters $\left[\frac{q_0+1}{2},\frac{q_0-3}{2},3\right]_{q_0}$, which is an MDS code. Moreover, $\rho(\widetilde{C_1(q_0)})=2$, so $\widetilde{C_1(q_0)}$ is both quasi-perfect and maximal. 

{\bf Case $s \ge 2$ even}: $\widetilde{C_s(q_0)}$ is a linear code with parameters $\left[\frac{q_0^s+1}{2},\frac{q_0^s+1}{2}-2s,4\right]_{q_0}$ with covering radius $\rho(\widetilde{C_s(q_0)})=3$, so $\widetilde{C_s(q_0)}$ is maximal.

{\bf Case $s \ge 3$ odd}: $\widetilde{C_s(q_0)}$ is a linear code with parameters $\left[\frac{q_0^s+1}{2},\frac{q_0^s+1}{2}-2s,3\right]_{q_0}$. The covering radius satisfies $\rho(\widetilde{C_s(q_0)})\in \{2,3\}$, so if $\rho(\widetilde{C_s(q_0)})=2$, then $\widetilde{C_s(q_0)}$ is quasi-perfect. In particular, when $s$ is odd and $3\leq s \leq s^*(q_0)$, then $d(\widetilde{C_s(q_0)})=3, \rho(\widetilde{C_s(q_0)})=2$, so $\widetilde{C_s(q_0)}$ is both quasi-perfect and maximal. 

\begin{theorem}
	Let $q_0 \ge 5$ be odd. 
	\begin{itemize}
		\item[(0)] If $s=1$, then $\widetilde{C_s(q_0)}$ is both quasi-perfect and maximal; 
		\item[(1)] If $s$ is odd and $3\leq s \leq s^*(q_0)$, then $\widetilde{C_s(q_0)}$ is both quasi-perfect and maximal; 
		\item[(2)] If $s \ge 2$ is even, then $\widetilde{C_s(q_0)}$ is maximal. 
	\end{itemize} 
\end{theorem}

The above results are summarized in the following table.

\begin{table}[h]
	\label{t2}
	\caption{ Parameters of $\widetilde{C_s(q_0)}$ in odd characteristic}
	\centering
	\begin{tabular}{|c|c|c|c|c|c|}
		\hline
		$q_0$ odd& $s$ & $d(\widetilde{C_s(q_0)})$ & $\rho(\widetilde{C_s(q_0)})$ & perfect/quasi-perfect & maximal \\ 
		\hline
		$q_0=3$ & $s \ge 2$ & $5$  & $3$ & quasi-perfect & $\checkmark$\\
		\hline
		$q_0 \ge 5$ & $1$ & $3$  & $2$ & quasi-perfect & $\checkmark$\\
\hline
		$q_0 \ge 5$& $s \geq 2$ even & $4$ & $3$ & - & $\checkmark$\\
		\hline
		$q_0 \ge 5$& $3\leq s \leq s^*(q_0)$, $s$ odd & $3$ & $2$ & quasi-perfect & $\checkmark$\\
		\hline
	\end{tabular}
\end{table}

\section{Concluding remarks}\label{sec8}

In this paper, we employed Weil-type estimates of character sums over finite fields to study the covering radius of generalized Zetterberg codes $C_s(q_0)$ for both even and odd $q_0$, and for each $q_0$, we determined the covering radius explicitly for a wide range of $s$. The study revealed many classes of quasi-perfect and maximal codes from $C_s(q_0)$. This complemented and strengthened previous works on odd characteristic. Here we raise one question which we might return in the future: 

Let $q_0$ be either even or odd and $s \ge 2$. It was proved in Corollary \ref{x3:cor} that $\rho(C_s(q_0)) \in \{2,3\}$. Moreover, by 
Propositions  \ref{x:prop-odd2}-\ref{x:prop-even} we have 
\begin{itemize}
    \item[(i)] if $s \ge 3$ is odd and $s \le s^*(q_0)$, then $\rho(C_s(q_0))=2$;

    \item[(ii)] if $s \ge 3$ is odd and $s \ge s_*(q_0)$, then $\rho(C_s(q_0))=3$. 
\end{itemize}
Here for odd and even $q_0$, the odd numbers $s^*(q_0)$ and $s_*(q_0)$ are defined explicitly at the end of Sections \ref{sec5} and \ref{sec6} respectively. 
 
{\bf Open question}: Can one determine additional values  of $\rho(C_s(q_0))$ for odd  $s$ in the range $\left(s^*(q_0),s_*(q_0)\right)$? Affirmative results would likely uncover new families of quasi-perfect and maximal codes arising from the generalized Zetterberg codes $C_s(q_0)$.

\bibliographystyle{IEEEtranS}

\bibliography{library(1)}

\end{document}